\documentclass{amsart}
\usepackage{amssymb,latexsym,graphicx,psfrag,epstopdf,hyperref,color}

\usepackage[all,arc,dvips]{xy}

\newcommand{\executeiffilenewer}[3]{%
\ifnum\pdfstrcmp{\pdffilemoddate{#1}}%
{\pdffilemoddate{#2}}>0%
{\immediate\write18{#3}}\fi%
}

\newtheorem{proposition}{Proposition}[section]
\newtheorem{theorem}[proposition]{Theorem}
\newtheorem{definition}[proposition]{Definition}

\newtheorem{lemma}[proposition]{Lemma}

\theoremstyle{remark}
\newtheorem{remark}{Remark}

\theoremstyle{remark}

\newtheorem*{acknowledgments}{Acknowledgments}

\def\frac#1#2{{#1\over#2}}

\numberwithin{equation}{section}

\def\co{\colon\thinspace}

\def\cal#1{\mathcal{#1}}

\begin{document}
\Large
\title{Small hyperbolic polyhedra}
\author[Shawn~Rafalski]{}
\address{Department of Mathematics and Computer Science, Fairfield University, Fairfield, CT 06824, USA}
\email{srafalski@fairfield.edu}
\date{\noindent January 2011. 
}

\maketitle

\centerline{\scshape shawn rafalski }
\medskip
{\footnotesize
\centerline{Department of Mathematics and Computer Science}
   \centerline{Fairfield University}
   \centerline{Fairfield, CT 06824, USA}
   \centerline{srafalski@fairfield.edu}
} 

\begin{abstract}
We classify the 3--dimensional hyperbolic polyhedral orbifolds that contain no embedded essential 2--suborbifolds, up to decomposition along embedded hyperbolic triangle orbifolds (turnovers). We give a necessary condition for a 3--dimensional hyperbolic polyhedral orbifold to contain an immersed (singular) hyperbolic turnover, we classify the triangle subgroups of the fundamental groups of orientable 3--dimensional hyperbolic tetrahedral orbifolds in the case when all of the vertices of the tetrahedra are non-finite, and we provide a conjectural classification of all the triangle subgroups of the fundamental groups of orientable 3--dimensional hyperbolic polyhedral orbifolds. Finally, we show that any triangle subgroup of a (non-orientable) 3--dimensional hyperbolic reflection group arises from a triangle reflection subgroup. \end{abstract}

\section{Introduction}\label{S:Intro}

Let $P$ be a finite volume 3--dimensional hyperbolic Coxeter polyhedron. That is, $P$ is the finite volume intersection of a finite collection of half-spaces in hyperbolic 3--space $\mathbb{H}^{3}$ in which the bounding planes of each pair of intersecting half-spaces meet at an angle of the form $\pi/n$, where $n\geq2$ is an integer (the geodesic of intersection is called an \emph{edge} of $P$, and the angle of intersection is called the \emph{dihedral angle} of $P$ along this edge). Then the group of isometries of $\mathbb{H}^{3}$ generated by the reflections in the faces of $P$ is a discrete group that acts on $\mathbb{H}^{3}$ with fundamental domain $P$. Let $\Gamma$ be the subgroup of index two in this reflection group generated by all the rotations of the form $r s$, where $r$ and $s$ are the reflections through two intersecting planes that support $P$.  We denote by $\cal{O}_{P}$ the quotient space $\mathbb{H}^{3}/\Gamma$. Then $\cal{O}_{P}$ is an orientable hyperbolic 3--orbifold called a \emph{hyperbolic polyhedral orbifold}. The group $\Gamma$ is sometimes denoted by $\pi_{1}(\cal{O}_{P})$ and called the \emph{fundamental group of} $\cal{O}_{P}$. We call $P$ a \emph{hyperbolic reflection polyhedron}.

\begin{figure}[htbp]
		\centering
		\def\svgwidth{1.5in}
	    	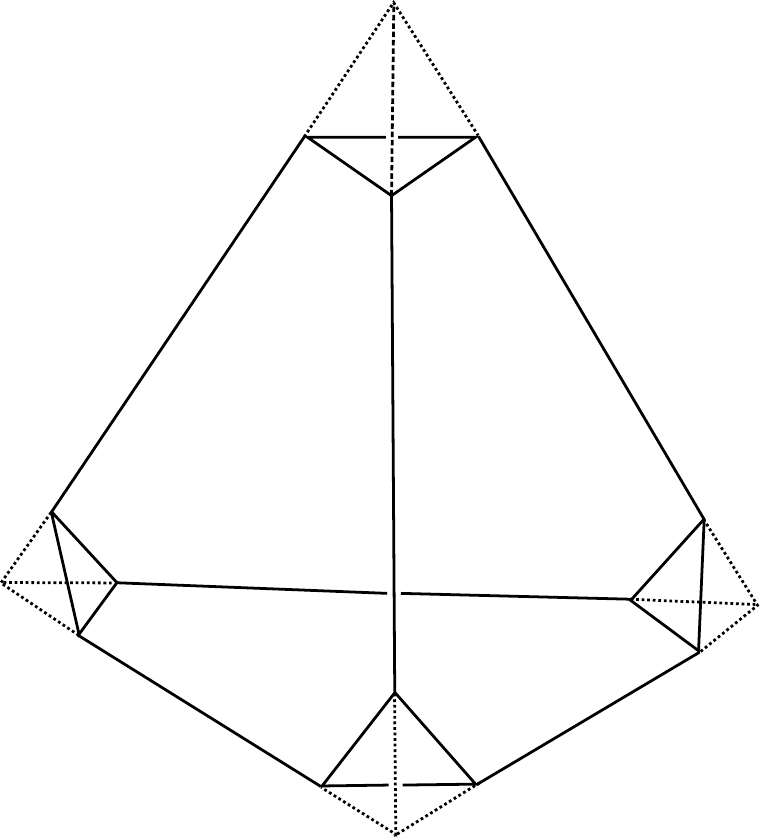
		\caption{The small Coxeter polyhedra in $\mathbb{H}^{3}$}
		\label{F:GenlzdTet}
\end{figure}

A \emph{small} hyperbolic reflection polyhedron corresponds to a hyperbolic 3--dimensional polyhedral orbifold that contains no embedded essential 2--suborbifolds, up to decomposition along embedded triangular 2--suborbifolds  (Definition \ref{D:smallpoly}). We classify these polyhedra (see Figure \ref{F:GenlzdTet}):

\begin{theorem}\label{T:SmallPolyhedra}
A 3--dimensional hyperbolic reflection polyhedron is small if and only if it is a generalized tetrahedron.
\end{theorem}

We also determine those hyperbolic polyhedral orbifolds that contain an immersed (singular) hyperbolic triangular 2--suborbifold. This result is a generalization of the partial classification of triangle groups inside of arithmetic hyperbolic tetrahedral reflection groups given by Maclachlan \cite{Maclachlan96-1}. In Section \ref{S:Turnovers}, we will provide a conjectural list of all the possibilities for immersed turnovers in all polyhedral orbifolds:

\begin{theorem}\label{T:TurnoversInPolyhedra}
If a hyperbolic polyhedral 3--orbifold contains a singular hyperbolic turnover that does not cover an embedded hyperbolic turnover, then at least one component of its Dunbar decomposition is a generalized tetrahedron, and the immersed turnover is contained in a unique such component. Furthermore, if $T$ is a generalized tetrahedron with all non-finite vertices and whose associated polyhedral 3--orbifold contains an immersed turnover, then, up to symmetry, $T$ is of the form $T[2,m,q;2,p,3]$ (in the notation described in Section \ref{S:Turnovers}) with $m \geq 6$, $q \geq 3$ and $p \geq 6$, and the immersed turnover has singular points of orders $m$, $q$ and $p$.
\end{theorem}

We also determine the triangle subgroups of 3--dimensional hyperbolic reflection groups as arising from triangle reflection subgroups:

\begin{theorem}\label{T:TriangleReflections}
Any (orientable) hyperbolic triangle subgroup of a (non-orientable) 3--dimensional hyperbolic reflection group $G$ arises as a subgroup of index two of a (non-orientable) hyperbolic triangle reflection subgroup of $G$. 
\end{theorem} 

Essential surfaces play an integral role in low-dimensional topology and geometry. One of the most important instances of this fact is the proof of Thurston's Hyperbolization Theorem for Haken 3--manifolds \cite{Thurston84-1}, \cite{Morgan84}. In brief, Thurston's Theorem is proved by decomposing a given 3--manifold $M$ (which is called \emph{Haken} if it contains an essential surface) along such surfaces as part of  a finite-step process that ends in topological solid balls, from which the hyperbolic structure on $M$ (whose existence is claimed by the theorem) is then, in a sense, reverse-engineered.

One difficulty that arises in attempting to extend the utility of essential surfaces to the orbifold setting is the possible presence of triangular hyperbolic 2--dimensional suborbifolds called \emph{hyperbolic turnovers}. For example, whereas an irreducible 3--manifold with non-empty and non-spherical boundary always contains an essential surface, this is not always the case in the orbifold setting, with hyperbolic turnovers presenting the principal barrier. Thurston's original definition of a Haken 3--orbifold was given for non-orientable 3--orbifolds with underlying space the 3--ball and with singular locus equal to the boundary of the ball \cite[Section 13.5, p. 324]{Thurstonnotes}. (The singular locus, in this instance, was meant to correspond to the boundary of a polyhedron.) Subsequent formulations of the definition of Haken (i.e., ``sufficiently large'' in \cite[Glossary]{Dunbar88-1} or ``Haken'' in \cite[Section 4.2, Remark]{BMP03-1}) were given for the orientable case and take into account the difficulties that arise from hyperbolic turnovers. Theorem \ref{T:SmallPolyhedra}, which is proved using the same observations that Thurston used to determine 3--orbifolds with the combinatorial type of a simplex as the original  ``non-Haken'' polyhedral orbifolds, echoes Thurston's original result \cite[Proposition 13.5.2]{Thurstonnotes}, with respect to this evolution of the language.

\begin{acknowledgments}
The author thanks Ian Agol for helpful conversations. Very special thanks to the referee for invaluable feedback and for recommending simplifications to some of the arguments.
\end{acknowledgments}


\section{Definitions}\label{S:Defns}

There are several excellent references for orbifolds \cite{BMP03-1}, \cite{CoopHodgKer00}. All of the 3--orbifolds considered in this paper are either orientable hyperbolic polyhedral 3--orbifolds or the result of cutting an orientable hyperbolic polyhedral 3--orbifold along a finite set of totally geodesic hyperbolic  turnovers or totally geodesic hyperbolic triangles with mirrored sides. A hyperbolic polyhedral 3--orbifold $\cal{O}_{P}$ is geometrically just two copies of its associated hyperbolic polyhedron $P$ with the corresponding sides of the two copies identified. Thus, $\cal{O}_{P}$ is a complete metric space of constant curvature $-1$ except along a 1--dimensional singular subset which is locally cone-like. If $P$ is compact, then  $\cal{O}_{P}$ is topologically a 3--sphere together with a trivalent planar graph  (corresponding to the 1--skeleton of $P$) with each edge marked by a positive integer to represent the submultiple of $\pi$  of  the dihedral angle at the corresponding edge of $P$. If $P$ is noncompact with finite volume, then its ideal vertices correspond to trivalent or quadrivalent vertices in the planar graph (again, corresponding to the 1--skeleton of $P$) and the sum of the reciprocals of the incident edge marks at each such vertex is equal to one or two, according to whether the vertex is trivalent or quadrivalent. In the noncompact case, $\cal{O}_{P}$ is topologically the result of taking a 3--sphere with this marked graph and removing a (closed)  3--ball neighborhood from each ideal vertex. The statements about the combinatorics of hyperbolic polyhedra in this paragraph are consequences of  Andreev's Theorem \cite{Andreev70-1}, \cite{Andreev70-2}, \cite{Roeder06}, \cite[Section 13.6]{Thurstonnotes}, \cite{Hodgson92-1}.

A (closed) orientable 2--orbifold is topologically a closed orientable surface with some finite set of its points marked by positive integers (greater than one). Every such 2--orbifold can be realized as a complete metric space of constant curvature with cone-like singularities at the marked points, and where the sign of the curvature depends only on the topology of the underlying surface together with the markings. A 2--orbifold is called \emph{spherical, Euclidean} or \emph{hyperbolic} according to the sign of its constant curvature realization. A \emph{turnover} is a 2--orbifold that is topologically a 2--sphere with three marked points, and a \emph{hyperbolic turnover} is a turnover for which the reciprocal sum of the integer markings is less than one. Although we will seldom deal with non-orientable objects, we define a \emph{hyperbolic triangle with mirrored sides} as a topological closed disk whose boundary is marked with three distinct points, each point labeled by an integer greater than one and such that the sum of the reciprocals of these integers is less than one, and with the connecting intervals in the boundary between these points marked as ``mirrors.'' Hyperbolic triangles with mirrored sides are non-orientable 2--orbifolds that are doubly covered by hyperbolic turnovers: they are the quotients of hyperbolic turnovers by an involution that fixes an embedded topological circle that passes through the marked points of the turnover.  Every embedded hyperbolic turnover in a hyperbolic 3--orbifold can either be made totally geodesic by an isotopy in the 3--orbifold (in which case the preimage in $\mathbb{H}^{3}$ under the covering map of this totally geodesic 2--suborbifold is a collection of disjoint planes, each tiled by a hyperbolic triangle that is determined by the markings of the singular points (e.g., \cite[Chapter IX.C]{Maskit87}, \cite[Theorem 2.1]{Adams05-1})) or else can be moved by an isotopy to be the boundary of a regular neighborhood of a totally geodesic hyperbolic triangle with mirrored sides.

An embedded orientable 2--suborbifold of $\cal{O}_{P}$ is topologically a surface that meets the marked graph transversely. We note that any simple closed curve $C \subset \partial P$ that meets the 1--skeleton transversely determines such a 2--suborbifold  by adjoining to $C$ the two topological disks that it bounds, one to either side of $\partial P \subset \cal{O}_{P}$. A closed path on $\partial P$ that is isotopic to a simple circuit in the dual graph to the 1--skeleton of $P$  is called a \emph{$k$--circuit}, where $k$ is the number of edges the path crosses. An embedded hyperbolic triangle with mirrored sides occurs as a suborbifold of $\cal{O}_{P}$ whenever $P$ has a triangular face all of whose edges are labeled two (in this case, the triangle with mirrored sides is topologically just the disc bounded by these three edges in the marked graph).

The terminology of this paragraph is introduced in terms of general orbifolds. A compact \emph{$n$--orbifold} $\cal{O}$ with boundary is a metrizable topological space which is locally diffeomorphic either to the quotient of $\mathbb{R}^{n}$ by a finite group action or to the quotient of $\mathbb{R}^{n-1} \times [0,\infty)$ by a finite group action, with points of the latter type making up the \emph{boundary} $\partial \cal{O}$ of $\cal{O}$ (itself an $(n-1)$--orbifold). We use the term \emph{orbifold ball} (respectively, \emph{orbifold disk}) to refer to the quotient of a compact 3--ball (respectively, 2--disk) by a finite group action.  We say a compact 3--orbifold $\cal{O}$ is \emph{irreducible} if every embedded spherical 2--suborbifold bounds an orbifold ball in $\cal{O}$. A 2--suborbifold $F \subset \cal{O}$ is called \emph{compressible} if either $F$ is spherical and bounds an orbifold ball or if there is a simple closed curve in $F$ that does not bound an orbifold disk in $F$ but that bounds an orbifold disk in $\cal{O}$, and \emph{incompressible} otherwise. There is a relative notion of $\partial$--incompressibility (whose exact definition we do not require). We call $F$ \emph{essential} if it is incompressible, $\partial$--incompressible and not parallel to a boundary component of $\cal{O}$. We call a compact irreducible 3--orbifold \emph{Haken} if it is either an orbifold ball, or a turnover crossed with an interval, or if it contains an essential 2--suborbifold but contains no essential turnover. A compact irreducible 3--orbifold is called \emph{small} if it contains no essential 2--suborbifolds and has (possibly empty) boundary consisting only of turnovers. (We note that a compact, orientable and irreducible orbifold is both Haken and small if and only if it is either a cone on a spherical turnover or a product of a turnover with an interval.) These definitions extend to any arbitrary 3--orbifold that is diffeomorphic to the interior of a compact 3--orbifold with boundary.  

We observe that Euclidean and hyperbolic turnovers are always incompressible because a simple closed curve on these objects always bounds an orbifold disk. As a consequence, in an irreducible 3--orbifold, any incompressible 2--orbifold (in fact, even any singular hyperbolic turnover) can be made disjoint from an embedded hyperbolic turnover. 

\begin{remark}\label{R:DunbarConvention} It is a consequence of a theorem proved by Dunbar that a hyperbolic polyhedral 3--orbifold can be decomposed (uniquely, up to isotopy) along a system of essential, pairwise non-parallel hyperbolic turnovers into pieces that contain no essential (embedded) turnovers, and, moreover, that each component of the decomposition is either a Haken or a small 3--orbifold  (\cite{Dunbar88-1}, \cite[Theorem 4.8]{BMP03-1}). An embedded hyperbolic turnover in a hyperbolic polyhedral 3--orbifold $\cal{O}_{P}$ will correspond to a simple closed curve in $\partial P$ that crosses exactly three edges whose dihedral angles sum to less than $\pi$. If such a curve is parallel in $\partial P$ to a triangular face of $P$ all of whose edges are labeled two, then the hyperbolic turnover corresponding to this curve is isotopic to the boundary of a regular neighborhood of a hyperbolic triangle with mirrored sides (the latter arising from the triangular face of $P$) in $\cal{O}_{P}$. In this case, one component of the Dunbar decomposition will consist of the regular neighborhood of this triangle with mirrored sides (in fact, this is a small 3--orbifold). The complement of this component in $\cal{O}_{P}$ is (orbifold) diffeomorphic to the complement of the triangle with mirrored sides in $\cal{O}_{P}$ (because the hyperbolic turnover collapses onto the mirrored triangle as the radius of the regular neighborhood goes to zero), and so, for convenience, we discard the component of the Dunbar decomposition corresponding to this regular neighborhood. 
\end{remark}

With the above convention in mind, we have the following:   

\begin{definition}\label{D:smallpoly}
A hyperbolic reflection polyhedron $P$ is \emph{small} if the Dunbar decomposition of $\cal{O}_{P}$ (with the convention of the preceding paragraph) consists of a single connected small component.
\end{definition}

In the projective model of $\mathbb{H}^{3}$, consider a linearly independent set of four points, any or all of which may lie on  the boundary of or outside of the projective ball. If the line segment between each pair of these points intersects the interior of the projective ball, then the points determine a \emph{generalized tetrahedron}. This polyhedron is obtained by taking the (possibly infinite volume) polyhedron in $\mathbb{H}^{3}$ spanned by the points and truncating its infinite volume ends by the dual hyperplanes to the super-ideal vertices. The resulting polyhedron has finite volume and all of its vertices are either finite or ideal. The faces arising from truncated super-ideal vertices---which are called, along with the finite and ideal vertices, \emph{generalized vertices}---are triangular, and the dihedral angle at each edge of these faces is $\pi/2$.  In particular, if a generalized tetrahedron $P$ is a Coxeter polyhedron, then any generalized vertex arising from a truncated face is a hyperbolic triangle that tiles (under the tiling associated to $P$)  a geodesic plane in $\mathbb{H}^{3}$ (and thus gives rise to an embedded hyperbolic triangle with mirrored sides in $\cal{O}_{P}$).

\section{Proof of Theorem \ref{T:SmallPolyhedra}}\label{S:MainTheoremProof}

\begin{figure}[htbp]
		\centering
		\def\svgwidth{1in}
	    	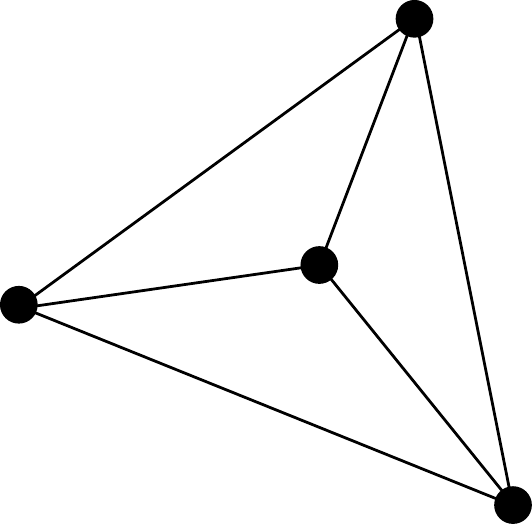
		\caption{The graph associated to a generalized tetrahedron}
		\label{F:GenTetGraph}
\end{figure}

\emph{Proof of \ref{T:SmallPolyhedra}.} Let $P$ be a 3--dimensional hyperbolic Coxeter polyhedron, and let $\cal{O}_{P}$ be its hyperbolic polyhedral 3--orbifold. First assume that $P$ is a generalized tetrahedron. Then $\cal{O}_{P}$ is topologically the 3--sphere with a marked planar graph as in Figure \ref{F:GenTetGraph}.    
Each dot in the figure represents a generalized vertex, and so is either a finite vertex, a triangle with mirrored sides or a Euclidean turnover cusp (the latter if the vertex is ideal). Any dot that represents a triangle corresponds to a non-separating hyperbolic turnover of the Dunbar decomposition of $\cal{O}_{P}$. Moreover, since any two hyperbolic turnovers can be made disjoint by an isotopy, any other turnovers in the Dunbar decomposition occur as topological 2--spheres that intersect the graph from the figure in exactly three distinct edges. But the only possibility for such a 2--sphere is one that surrounds a dot, and that therefore is parallel to a generalized vertex of $P$. So the Dunbar decomposition of $\cal{O}_{P}$ (under the convention of Definition \ref{D:smallpoly}) has a single component. 

To see that this component is small, we consider the graph of Figure \ref{F:GenTetGraph} as the 1--skeleton of a tetrahedron in the 3--sphere. Using standard topology arguments, it can be shown that an incompressible 2--suborbifold intersects the interior of this tetrahedron in triangles and quadrilaterals. But a triangular intersection  implies that the incompressible 2--suborbifold is isotopic to the hyperbolic turnover associated to a generalized vertex, and a quadrilateral intersection produces a compression. So $P$ is small if it is a generalized tetrahedron.

Now assume that $P$ is  small. The rest of the proof of Theorem \ref{T:SmallPolyhedra} depends on the following simple observation \cite[Proposition 13.5.2]{Thurstonnotes}:

\begin{remark}\label{R:Inc2Orbs}
Suppose that $C \subset \partial P$ is a simple closed curve that is transverse to, forms no bigons with, does not surround a single vertex of, and that crosses at least two distinct edges of  the 1--skeleton of $P$. Then $C$ determines an incompressible 2--suborbifold of $\cal{O}_{P}$ if and only if (1) it intersects any face in a connected set or not at all and (2) it intersects the common edge of two adjacent faces whenever its intersection with both faces is nonempty.
\end{remark}

We begin with the following fact about triangular faces of $P$:

\begin{lemma}\label{L:NoPrisms}
If $T$ is a triangular face of $P$, then $T$ corresponds to a hyperbolic turnover in $\cal{O}_{P}$ or $P$ is a generalized tetrahedron.
\end{lemma}

\begin{figure}[htbp]
		\centering
		\def\svgwidth{4.5in}
	    	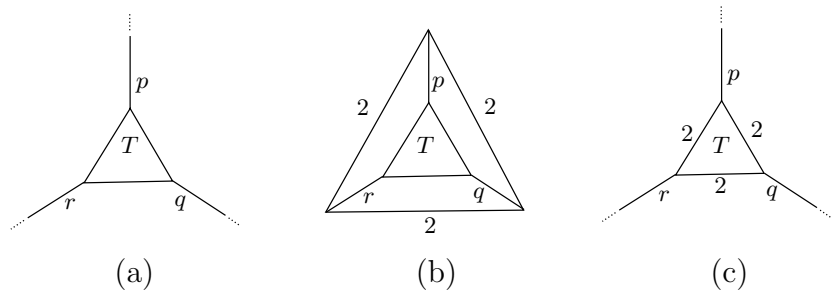
		\caption{Triangular faces}
		\label{F:NoPrisms}
\end{figure}

\emph{Proof of \ref{L:NoPrisms}.}
Suppose that $T$ is as in Figure \ref{F:NoPrisms}a (in this and all subsequent figures in this section, we depict $P$ by a planar projection). If $1/p + 1/q + 1/r  \geq 1$, then the three edges incident to the vertices of $T$ must intersect (or meet at a Euclidean turnover) \cite[Lemmata 3.2 and 3.3]{Roeder06}, in which case $P$ is a generalized tetrahedron (possibly with an ideal vertex). Otherwise, we have $1/p+1/q+1/r < 1$. Then the 3--circuit around this face determines a hyperbolic turnover in $\cal{O}_{P}$ whose associated triangle in $P$ must be boundary-parallel (in $P$) because $P$ is small. The two possibilities are shown in Figures \ref{F:NoPrisms}b (in which the hyperbolic turnover collapses to the outermost face) and \ref{F:NoPrisms}c (in which the hyperbolic turnover collapses to $T$).  \hfill  \fbox{\ref{L:NoPrisms}}\\

Throughout the rest of the proof, we will use the observation from the above lemma, i.e., that any 3--circuit in a small hyperbolic polyhedron surrounds a generalized vertex. In the case when the 3--circuit determines a hyperbolic turnover, this follows by the fact that a hyperbolic turnover in a hyperbolic 3--orbifold always corresponds to a totally geodesic 2--suborbifold (according to the second paragraph in Section \ref{S:Defns}; compare also with the incompressibility observation in the paragraph preceding Remark \ref{R:DunbarConvention}): Because the polyhedron $P$ is small, this totally geodesic 2--suborbifold of $\cal{O}_{P}$ cannot be an embedded hyperbolic turnover (because $\cal{O}_{P}$ has no boundary, and so such a turnover would have to be essential), and therefore must be a triangle with mirrored sides that corresponds to a triangular face of $P$.
 
\begin{figure}[htbp]
		\centering
		\def\svgwidth{4in}
	    	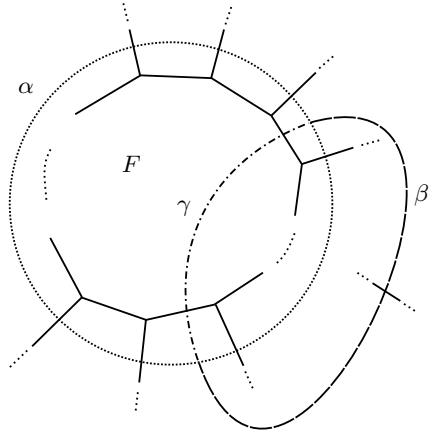
		\caption{A face of $P$ and a compression}
		\label{F:GeneralNFace1}
\end{figure}

\begin{figure}[htbp]
		\centering
		\def\svgwidth{6in}
	    	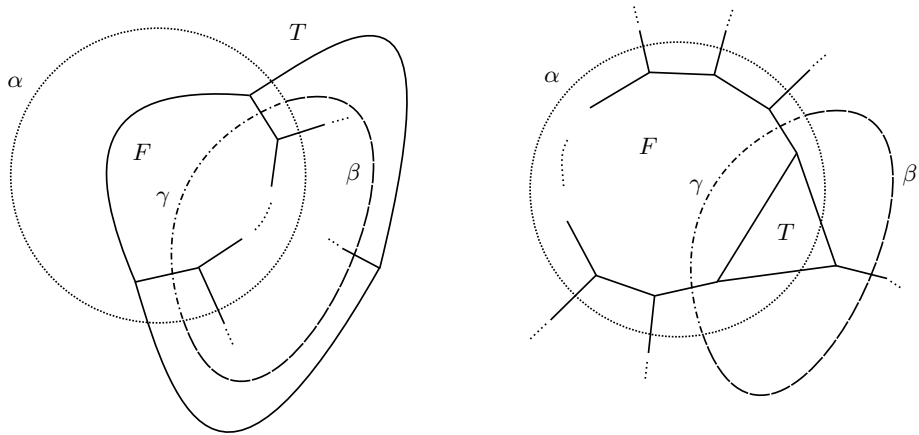
		\caption{Two projections of a  face of $P$ with adjacent triangle}
		\label{F:GeneralNFace2}
\end{figure}

Consider an $n$--sided face $F$ of $P$, as in Figure \ref{F:GeneralNFace1}. 
Assume that $n\geq4$. The $n$--circuit $\alpha$ around $F$ determines a 2--orbifold that must be compressible, with a compressing orbifold disk whose intersection with $\partial P$ appears as the dashed arc $\beta$ in the figure. Since $n \geq 4$, it must be that each side of the 3--circuit  $\beta \cup \gamma$ contains at least two edges radiating outward from $F$ (that is, edges meeting $F$ only in vertices). Since $\cal{O}_{P}$ is small, $\beta \cup \gamma$ bounds a triangle $T \subset \partial P$. Figure \ref{F:GeneralNFace2} illustrates the two possibilities, depending on the side of $\beta \cup \gamma$ to which $T$ lies. Of course, these differ only by the choice of projection of $P$ into the plane. 

\begin{figure}[htbp]
		\centering
		\def\svgwidth{2in}
	    	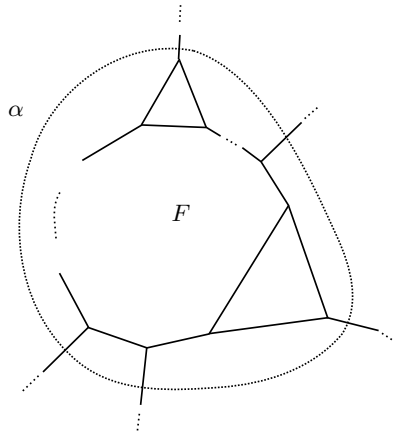
		\caption{A face of $P$ with all of its adjacent triangles, and a $k$--circuit}
		\label{F:GeneralNFaceWithTriangles}
\end{figure}

We now consider all such compressions of this 2--orbifold, and  all of the resulting adjacent triangles to $F$. Let $\alpha$ denote the $k$--circuit that encloses $F$ and these triangles, as in Figure \ref{F:GeneralNFaceWithTriangles}.

\begin{figure}[htbp]
		\centering
		\def\svgwidth{2.5in}
	    	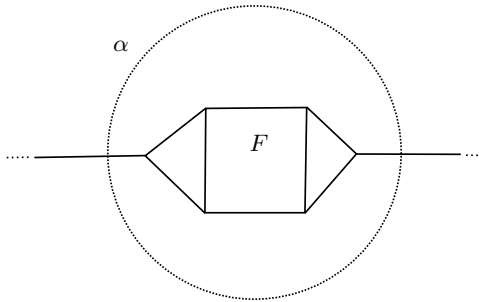
		\caption{A face of $P$ with all of its adjacent triangles, and a 2--circuit}
		\label{F:QuadFace1}
\end{figure}

If $k =2$, then $F$ must be a quadrilateral with two triangles adjacent to it on opposite sides, in which case $P$ is a triangular prism with one face that corresponds to a hyperbolic turnover in $\cal{O}_{P}$ as in Lemma \ref{L:NoPrisms}, i.e., $P$ is a generalized tetrahedron. See Figure \ref{F:QuadFace1}. 

\begin{figure}[htbp]
		\centering
		\def\svgwidth{1.7in}
	    	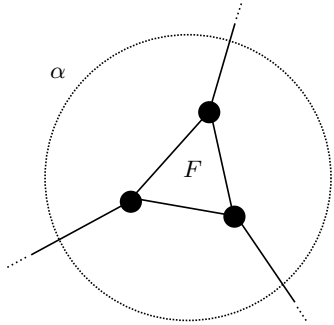
		\caption{A face of $P$ with all of its adjacent triangles, and a 3--circuit}
		\label{F:GeneralNFaceEndGame}
\end{figure}

If $k=3$, then $\alpha$ surrounds a generalized vertex to the outside. In this case, the face $F$ must be as in Figure \ref{F:GeneralNFaceEndGame}, where each dot represents either a finite vertex, an ideal vertex or a hyperbolic triangle. Filling in the generalized vertex to the outside of $\alpha$, we have that $P$ is a generalized tetrahedron.

If $k>3$, then the 2--orbifold determined by $\alpha$ has a compression. But any such compression would add an adjacent triangular face to $F$, and we have assumed that $\alpha$ encloses all such triangles. So $k \leq 3$. This completes the proof.  \hfill  \fbox{\ref{T:SmallPolyhedra}}\\

\section{Turnovers in Hyperbolic Polyhedra}\label{S:Turnovers}

In this final section, we prove Theorems \ref{T:TurnoversInPolyhedra} and \ref{T:TriangleReflections}, and provide a classification of the immersed hyperbolic turnovers in those tetrahedral orbifolds that arise from tetrahedra with no finite vertices. Although Theorem \ref{T:TriangleReflections} does not follow from Theorem \ref{T:TurnoversInPolyhedra}, we will provide the proof of the former in the midst of the proof of the latter, as it contains an observation that is necessary for both proofs.

The author showed that if a hyperbolic 3--orbifold contains a  singular hyperbolic turnover, then that turnover must be contained in a low-volume small 3--suborbifold \cite{Rafalski10}. In particular, we have the following

\begin{theorem}\label{T:JSJDecompEsque}
(Theorem 1.1 and Corollary 1.3 in \cite{Rafalski10}) Let $Q$ be a compact, irreducible, orientable, atoroidal 3--orbifold. Then any immersion $f\co \cal{T} \to Q$ of a hyperbolic turnover into $Q$ is homotopic into a unique component of the Dunbar decomposition of $Q$, up to covers of parallel boundary components of the decomposition. Moreover, if $f$ is a singular immersion that does not cover an embedded turnover or triangle with mirrored sides, then the component containing $f(\cal{T})$ is unique, and it is a small 3--orbifold.
\end{theorem}

\emph{Proof of \ref{T:TurnoversInPolyhedra}.}
If $\cal{O}_{P}$ is a hyperbolic polyhedral 3--orbifold, then it is homeomorphic to the interior of an orbifold that satisfies the hypotheses of Theorem \ref{T:JSJDecompEsque}. If $\cal{O}_{P}$ contains a singular turnover, then this turnover is contained in a small component of the Dunbar decomposition of $\cal{O}_{P}$, and Theorem \ref{T:SmallPolyhedra} classifies these small orbifolds as generalized tetrahedral orbifolds. 

It remains to provide a classification of the generalized tetrahedra whose associated 3--orbifolds contain immersed turnovers. We will do so for generalized tetrahedra all of whose vertices are non-finite. \emph{See the summary at the end of the paper for the results of the classification.} The techniques we use to provide this classification can be used to classify the immersed turnovers in all tetrahedral orbifolds, thereby extending and completing the classification begun by Maclachlan in the case of compact (non-generalized) tetrahedral orbifolds \cite{Maclachlan96-1}, however, the case-by-case analysis required to complete this classification in general is somewhat excessive.

\begin{figure}[htbp]
		\centering
		\def\svgwidth{3in}
	    	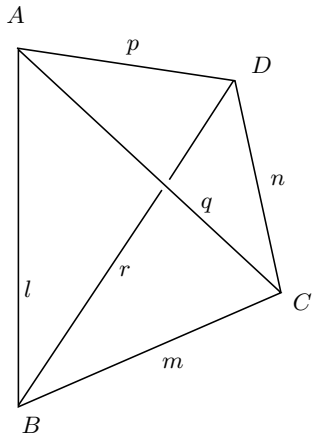
		\caption{$T[l,m,q;n,p,r]$}
		\label{F:TLMQNPR}
\end{figure}

We let $T[l,m,q;n,p,r]$ denote the hyperbolic generalized tetrahedron $ABCD$ with dihedral angles $\pi/l$, $\pi/m$, $\pi/q$, $\pi/n$, $\pi/p$ and $\pi/r$, as in Figure \ref{F:TLMQNPR}, with the convention that a vertex of $T$ is truncated (respectively, ideal) if the dihedral angles of its three coincident edges sum to less than (respectively, equal to) $\pi$.
The conditions on $l, m, q, n, p$ and $r$ guaranteeing the existence of $T[l,m,q;n,p,r]$ are known \cite{Ushijima06}. In particular, there are nine compact (non-truncated) tetrahedra  (e.g. \cite[Chapter 7]{Ratcliffe94-1}),  all of whose associated orbifolds contain singular turnovers. We note that, of the nine compact (non-truncated) tetrahedra, eight yield arithmetic hyperbolic 3--orbifolds. As we noted above, Maclachlan classified almost all of the immersed turnovers in these arithmetic tetrahedral orbifolds using arithmetic methods. Our geometric technique can be considered as an alternative means to prove (and extend) those results, without appeal to arithmeticity.

Denote by $\cal{O}_{T}$ the 3--orbifold determined by $T[l,m,q;n,p,r]$. Recall from Section \ref{S:Defns} that any hyperbolic turnover in a hyperbolic 3--orbifold that does not collapse onto a hyperbolic triangle with mirrored sides may be assumed to be totally geodesic. It also follows from the incompressibility of hyperbolic turnovers in irreducible orbifolds that an immersed turnover must be disjoint from any embedded turnover \cite[Lemma 5.3]{Rafalski10}. Consequently, if $\cal{T}$ is a hyperbolic turnover, then an immersion $f\co \cal{T} \to \cal{O}_{T}$ lifts to the universal cover $\mathbb{H}^{3}$ as a collection of geodesic planes with some intersections---two or more of these planes will intersect whenever there is a covering transformation (i.e., an element of the fundamental group $\pi_{1}(\cal{O}_{T})$ of $\cal{O}_{T}$, which is just the group of isometries of $\mathbb{H}^{3}$ that yields the quotient $\cal{O}_{T}$) that does not move one plane completely disjoint from some of the others, and this must occur if there is a singular immersion of a turnover in $\cal{O}_{T}$---and, additionally, the collection of planes determined by an immersed turnover must be disjoint from the collection of planes determined by any turnover corresponding to a generalized vertex of $T$.

After these preliminaries, we now will prove Theorem \ref{T:TriangleReflections}.

\emph{Proof of \ref{T:TriangleReflections}.} Let $P \subset \mathbb{H}^{3}$ be a polyhedron that generates the non-orientable 3--dimensional hyperbolic polyhedral reflection group $G$, and let $S \subset G$ be an orientable triangle subgroup. Then $S$ is generated by two elliptic elements in $G$ and  stabilizes a plane $\Pi_{S} \subset \mathbb{H}^{3}$. In particular, $\Pi_{S}$ meets the axis of every element of $S$ at a right angle, and the intersections of $\Pi_{S}$ with these axes comprise the vertex set of a tiling of $\Pi_{S}$ by hyperbolic triangles. Every such vertex will have $k$ lines passing through it (where $k$ is the order of the elliptic element stabilizing the vertex) that are the perpendicular intersections with $\Pi_{S}$ of $G$--translates of a face of $P$. This set of lines and their intersections generates a tiling of $\Pi_{S}$ by hyperbolic triangles that corresponds to a hyperbolic triangle with mirrored sides in the non-orientable hyperbolic orbifold $\mathbb{H}^{3}/G$,  and this 2--orbifold is covered by the hyperbolic turnover corresponding to $S$. Therefore, $S$ is contained in the triangle reflection subgroup of $G$ that corresponds to this non-orientable triangle 2--orbifold.   \hfill \fbox{\ref{T:TriangleReflections}}\\

We take a moment to emphasize the observation from the above proof: Any maximal (orientable) triangle subgroup of 3--dimensional hyperbolic polyhedral reflection group has as a fundamental domain a triangle whose edges are contained in the faces of the corresponding polyhedral tiling of $\mathbb{H}^{3}$ (the edges may intersect multiple faces of the polyhedral tiling). This fact is used in the next paragraph.

\emph{Here is the strategy for classifying the immersed turnovers of $\cal{O}_{T}$ (the proof is somewhat lengthy, but this paragraph contains the core idea):} Let $\cal{T}$ be a hyperbolic turnover. Up to conjugacy, there is a unique discrete orientation-preserving group of isometries of the hyperbolic plane $\mathbb{H}^{2}$ corresponding to the tiling of $\mathbb{H}^{2}$ by copies of the triangle that determines $\cal{T}$ (the fundamental group $\pi_{1}(\cal{T})$ of $\cal{T}$). If $f \co \cal{T} \to \cal{O}_{T}$ is an immersion, then $f$ may be assumed to have totally geodesic image. Consider a plane $\Pi_{\cal{T}}$ in the collection of planes in $\mathbb{H}^{3}$ corresponding to $f(\cal{T})$. This plane is stabilized by a copy of the fundamental group of some turnover (possibly a smaller turnover that is covered by $f(\cal{T})$, if the fundamental group of $f(\cal{T})$ is not maximal)---a subgroup $\Gamma$ of the fundamental group of the orbifold $\cal{O}_{T}$---for which there is a tiling of $\Pi_{\cal{T}}$ by hyperbolic triangles whose edges are a (possibly proper) subset of the intersections of $\Pi_{\cal{T}}$ with $\Gamma$--translates of the faces of $T$, and whose vertices are a (possibly proper) subset of the perpendicular intersections of $\Pi_{\cal{T}}$ with $\Gamma$--translates of the edges of $T$. We will locate all of the immersed turnovers in $\cal{O}_{T}$ by reversing this process, that is, by determining exactly the hyperbolic planes in the universal cover $\mathbb{H}^{3}$ that are stabilized by a triangle subgroup of $\pi_{1}(\cal{O}_{T})$.  We therefore choose an arbitrary edge $e_{1}$ of $T$ and develop copies of $T$ in $\mathbb{H}^{3}$ (by reflecting in faces) until we find another edge $e_{2}$ which is coplanar with but which shares no (generalized) vertex with $e_{1}$. Since we need only concern ourselves with maximal triangle subgroups, the observation following the proof of Theorem \ref{T:TriangleReflections} allows to assume that the common plane, which we denote by $\Pi_{F}$ (where $F$ is a face of $T$ incident to $e_{1}$), consist of developed faces of $T$.  Let $\Pi_{1}$ be the plane containing another face of $T$ incident with $e_{1}$, and let $\Pi_{2}$ be the plane containing another face of (a developed image of) $T$ containing $e_{2}$. Suppose that $\Pi_{1}$ and $\Pi_{2}$ intersect $\Pi_{F}$ at angles of $\pi/a$ and $\pi/b$, respectively.  If $\Pi_{1}$ and $\Pi_{2}$ intersect at an angle of $\pi/c$, and if $1/a+1/b+1/c < 1$, then the rotations about edges $e_{1}$ and $e_{2}$ (of orders $a$ and $b$, respectively), will generate a triangle subgroup of $\pi_{1}(\cal{O}_{T})$, and the invariant plane for that subgroup will project to an immersed turnover in $\cal{O}_{T}$ (every developed edge of $T$ that intersects the invariant plane for this triangle subgroup at an oblique angle will correspond to an immersion of the turnover). This determines a maximal triangle subgroup of $\pi_{1}(\cal{O}_{T})$, and the type of the corresponding immersed turnover will be $(a,b,c)$. In most cases, we will show that there can be no such edge $e_2$ that is both coplanar with $e_1$ and that has an incident face whose corresponding plane $\Pi_2$ intersects the plane $\Pi_1$, which rules out the possibility of an immersed turnover. In the other cases, we will find a turnover after a minimal development of $T$. Thus, our determination of the immersed turnovers in $\cal{O}_{T}$ will be complete.

We divide the remainder of the proof of Theorem \ref{T:TurnoversInPolyhedra} into subsections.

\subsection{The case when a single edge separates $e_{1}$ from $e_{2}$:}\label{SS:SingleEdgeCrossed} 

\subsubsection{The single separating edge has order 2:}\label{SSS:SingleEdgeCrossedOrder2} 
\begin{figure}
	\centering
	\def\svgwidth{3.5in}
	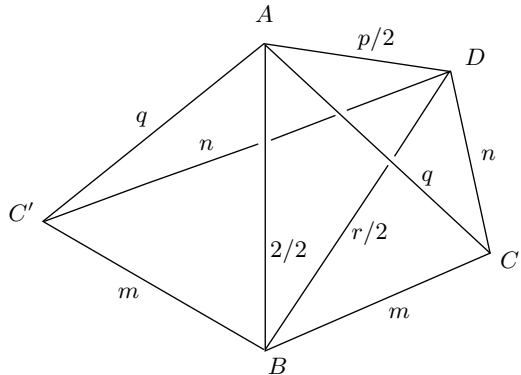
	\caption{Two copies of the tetrahedron $T[2,m,q;n,p,r]$}
	\label{F:Tet1EdgeOrder2}
\end{figure}

To begin, we determine the case in which the immersed turnover can be found after crossing only one edge between $e_{1}$ and $e_{2}$ (there must be at least one edge crossed, in this process, to ensure that the turnover is not parallel to a vertex). Consider Figure \ref{F:Tet1EdgeOrder2}, which shows two copies of the tetrahedron $T[2,m,q;n,p,r]$. Each edge is labeled according to the submultiple of $\pi$ for the dihedral angle there (so, for example, the edge $AD$ has a dihedral angle of $2\pi/p$). In particular, the points $A$, $B$, $C$ and $C'$ are coplanar. We use $F$ to denote the face $ABC$ of $T$ and $\Pi_{F}$ to denote the plane that contains $F$. We consider the edges $e_{1}=AC'$ and $e_{2}=BC$, and the planes $\Pi_{1}=AC'D$, $\Pi_{F}$ and $\Pi_{2}=BCD$. Under the assumption that all of the vertices of $T$ are non-finite, we observe it is necessary for $m, q, p$ and $r$ to all be at least 3. From the figure, we see that $\Pi_{1}$ meets $\Pi_{F}$ at an angle of $\pi/q$ and that $\Pi_{F}$ meets $\Pi_{2}$ at an angle of $\pi/m$, and so we are left to determine whether or not $\Pi_{1}$ and $\Pi_{2}$ intersect, and at what angle this possible intersection occurs.
The vertex $D$ is either ideal or truncated. If it is ideal, then its link is the orbifold quotient of a horosphere by a Euclidean triangle group. If it is truncated, then it corresponds to a geodesic plane in the universal cover that is stabilized by a hyperbolic triangle group. In both cases, we illustrate the straightforward geometric determination of the conditions on $n, p$ and $r$ that ensure the intersection of $\Pi_{1}$ and $\Pi_{2}$ in the link of $D$, and determine the angle at which any intersection occurs \cite[Section 9.4]{Rafalski10}.

\begin{figure}
	\centering
	\def\svgwidth{5.5in}
	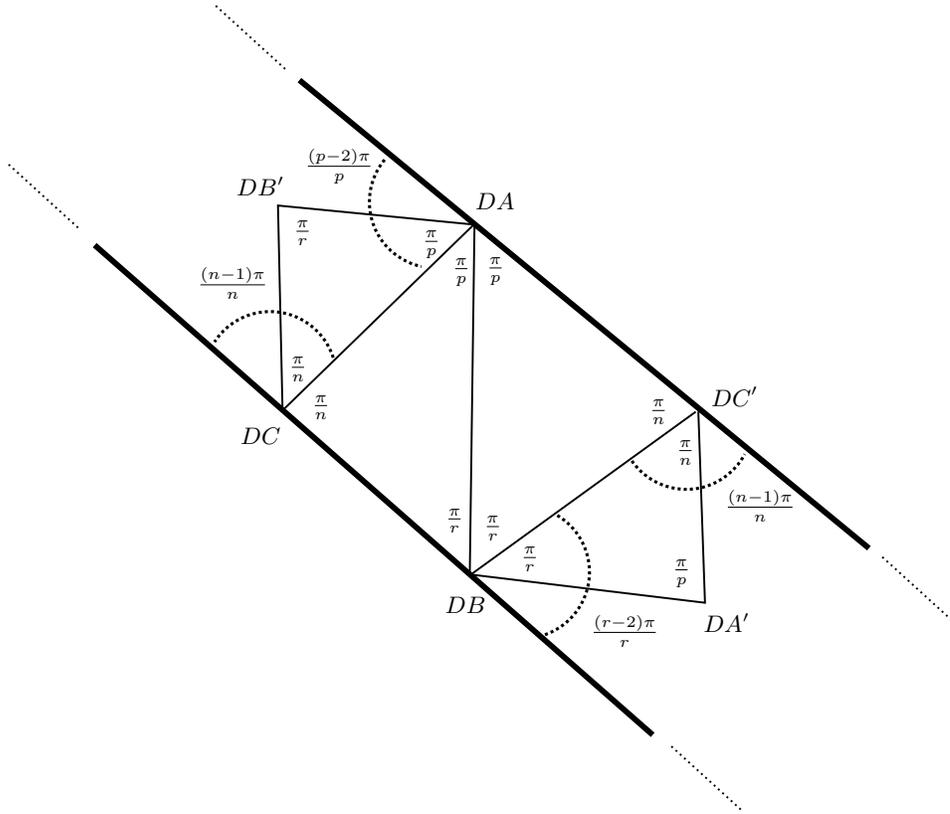
		\caption{Part of the link of a non-finite vertex of $T[2,m,q;n,p,r]$}
		\label{F:Tet1EdgeOrder2Link}
\end{figure}

Consider Figure \ref{F:Tet1EdgeOrder2Link}, which illustrates part of the link of $D$ \emph{as viewed from} $D$ (this is either a hyperbolic plane or a Euclidean plane corresponding to the horosphere centered at an ideal vertex). 
The vertices in the picture are labeled according to the edges of $T$ that are incident at $D$ (the labels $DA'$ and $DB'$ represent edges in the development of $T$ that are the reflections of the edges $DA$ and $DB$ through the faces $BC'D$ and $ACD$, respectively, in Figure \ref{F:Tet1EdgeOrder2}). Assume first that $n >2$. Because $p$ must be at least 3 (and similarly for $r$), we have the inequality $(p-2)\pi/p+(n-1)\pi/n \geq \pi$ (and similarly $(r-2)\pi/r+(n-1)\pi/n \geq \pi$). The angles with the measures from the previous sentence are indicated in the figure as the labels of the four \emph{dotted} arcs (all other angles in the figure refer to the measure at the appropriate triangular vertex). Using this inequality, we conclude that the indicated bold rays directed northwest from $DA$ and $DC$ do not intersect, because the sum of the angles that these rays make with the segment from $DC$ to $DA$ is at least $\pi$ (and similarly for the rays directed southeast from $DC'$ and $DB$, because the sum of the angles that these rays make with the segment from $DB$ to $DC'$ is at least $\pi$). Consequently, the bold lines in the figure (and the corresponding planes $\Pi_{1}$ and $\Pi_{2}$) cannot intersect in this case. A similar argument implying that $\Pi_{1}$ and $\Pi_{2}$ do not intersect holds when $n=2$ and both $p$ and $r$ are greater than 3: The rays directed northwest from $DA$ and $DC$ make angles with the segment between these two points of $(p-2)\pi/p \geq \pi/2$ and $\pi/2$, respectively, and so the sum of these angles will be at least $\pi$ (when $n=2$ and $r\geq4$, the same argument proves that the southeast rays from $DC'$ and $DB$ do not intersect). Finally, if $n=2$ and $p=3$ (respectively, $r=3$), then it is easily seen $\Pi_{1}$ and $\Pi_{2}$ intersect at an angle of $\pi/r$ (respectively, $\pi/p$), and the line of intersection passes through the point $DB'$ (respectively, $DA'$).

\begin{figure}
	\centering
	\def\svgwidth{4.5in}
	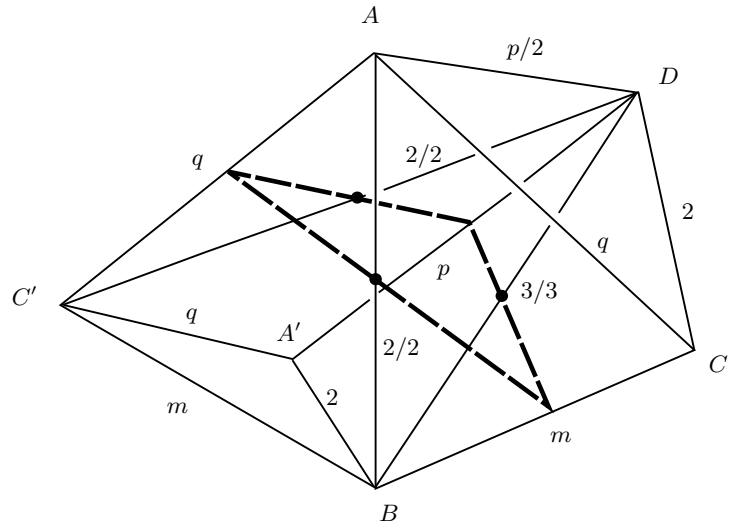
		\caption{A $(q,m,p)$ triangle in $T[2,m,q;2,p,3]$}
		\label{F:Tet1NoSupergroup}
\end{figure}

We therefore have, when $l=2$ and our search for a turnover crosses only one edge, that an immersed turnover only arises when $n=2$ and either $r=3$ or $p=3$. If $r=3$, then this yields a triple of planes intersecting pairwise in angles of $\pi/q$, $\pi/m$ and $\pi/p$, with $q \geq 3$, $m \geq 6$ and $p \geq 6$. If $p=3$, then the pairwise angles of intersection are $\pi/q$, $\pi/m$ and $\pi/r$, with $q \geq 6$, $m \geq 3$ and $r \geq 6$. (The inequalities are induced by the assumption that all of the vertices of $T$ are non-finite.) By analyzing Table \ref{Ta:TurnoverSupergroups} (whose data is collected from Singerman \cite{Singerman72}), we see that this triple of planes does not yield a triangle group that contains any other triangle group. By comparing the second column of the table with the first, we note that it is possible for this triple of planes to yield a triangle group that is contained in some larger triangle group. However, it is not possible for such a supergroup to be a subgroup of $\pi_{1}(\cal{O}_{T})$. This follows from the observation in the paragraph following the proof of Theorem \ref{T:TriangleReflections}: Because such a supergroup would be a maximal triangle subgroup of $\pi_{1}(\cal{O}_{T})$ stabilizing the plane that contains the $(q,m,p)$ (or $(q,m,r)$) triangle, there would have to be edges in the development of $T$ that intersect the interior of the $(q,m,p)$ (or $(q,m,r)$) triangle perpendicularly (these intersections would be necessary for the corresponding orbifold covering of the smaller turnover by the larger $(q,m,p)$ or $(q,m,r)$ turnover). By  construction, there are no such perpendicular intersections in the interior of the triangle. See Figure \ref{F:Tet1NoSupergroup}, which illustrates the case when $r=3$. As can be seen in the figure, no developed edges of $T$ intersect the interior of the $(q,m,p)$ triangle (the intersections with this triangle that yield immersions of the corresponding turnover are indicated by the dots). Consequently, we can conclude that the $(q,m,p)$ or $(q,m,r)$ triangle determined by $\Pi_{1}, \Pi_{F}$ and $\Pi_{2}$ is not parallel to any of the vertices of $T$, and therefore that it determines an immersed turnover in $\cal{O}_{T}$, because $\cal{O}_{T}$ is small. The observations of this paragraph are summarized in items (1) and (2) at the conclusion of the paper.
\begin{table}
  \begin{center}
    \begin{tabular}{ | c | c | c | c | }
      \hline
      \rule[-8pt]{0pt}{22pt}
      Supergroup	&	Subgroup	&	Index	&	Normal	\\ \hline
      $(3,3,t)$	&	$(t,t,t)$		&	$3$		&	Yes		\\ \hline
      $(2,3,2t)$	&	$(t,t,t)$		&	$6$		&	Yes		\\ \hline
      $(2,s,2t)$	&	$(s,s,t)$		&	$2$		&	Yes		\\ \hline
      $(2,3,7)$	&	$(7,7,7)$		&	$24$		&	No		\\ \hline
      $(2,3,7)$	&	$(2,7,7)$		&	$9$		&	No		\\ \hline
      $(2,3,7)$	&	$(3,3,7)$		&	$8$		&	No		\\ \hline
      $(2,3,8)$	&	$(4,8,8)$		&	$12$		&	No		\\ \hline
      $(2,3,8)$	&	$(3,8,8)$		&	$10$		&	No		\\ \hline
      $(2,3,9)$	&	$(9,9,9)$		&	$12$		&	No		\\ \hline
      $(2,4,5)$	&	$(4,4,5)$		&	$6$		&	No		\\ \hline
      $(2,3,4t)$	&	$(t,4t,4t)$		&	$6$		&	No		\\ \hline
      $(2,4,2t)$	&	$(t,2t,2t)$		&	$4$		&	No		\\ \hline
      $(2,3,3t)$	&	$(3,t,3t)$		&	$4$		&	No		\\ \hline
      $(2,3,2t)$	&	$(2,t,2t)$		&	$3$		&	No		\\ \hline
    \end{tabular}
    \medskip
    \caption{Triangle Supergroups and Subgroups}\label{Ta:TurnoverSupergroups}
  \end{center}
\end{table}

\subsubsection{The single separating edge has order 3:}\label{SSS:SingleEdgeCrossedOrder3} 
\begin{figure}
	\centering
	\def\svgwidth{4.5in}
	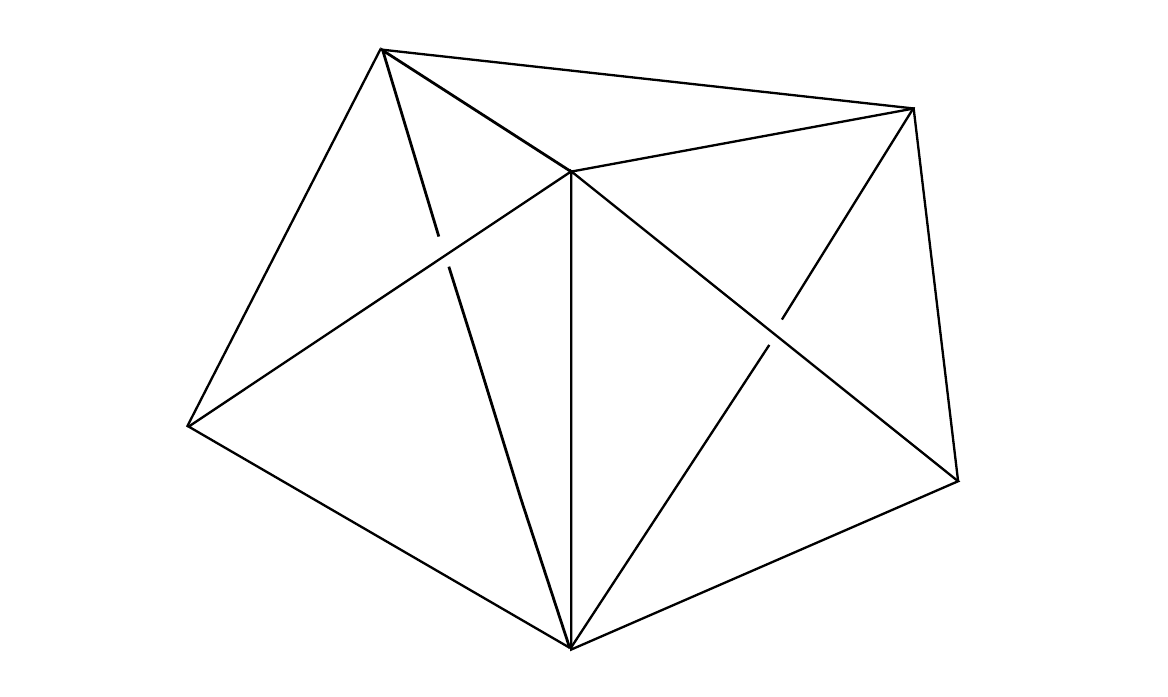
	\caption{Three copies of the tetrahedron $T[3,m,q;n,p,r]$}
	\label{F:Tet1EdgeOrder3}
\end{figure}

We next turn to the case in which the immersed turnover can be found after crossing only one edge between $e_{1}$ and $e_{2}$, where the order of the crossed edge is $l=3$. See Figure \ref{F:Tet1EdgeOrder3}.
Let $e_{1}=AD'$, $e_{2}=BC$, $\Pi_{1}=AC'D'$ and $\Pi_{2}=BCD$. We make several preliminary observations:
\begin{enumerate}
	\item Any two (distinct) planes that truncate developed vertices will always be disjoint.
	\item By (1) and by the fact that $T$ has no finite vertices, any two developed edges of the tetrahedron (whose corresponding geodesics in $\mathbb{H}^{3}$ are distinct) will always be disjoint.  
	\item It is always the case that the plane containing one face of a generalized hyperbolic tetrahedron will be disjoint from the plane that truncates the vertex opposite to that face. 
	\item By (2), if two planes corresponding to two developed faces of $T$ meet a third plane that corresponds to a developed face of $T$, then any intersection of the first two planes must occur on the side of the third plane where the two interior supplementary angles of intersection sum to less than $\pi$. 
	\item Any two planes corresponding to two developed faces that both intersect a third plane that truncates a developed vertex intersect if and only if their intersections with that truncated plane (i.e., with the link of the generalized vertex) do so. A corresponding statement is also true in the case when the developed vertex is ideal, that is, that two planes corresponding to two developed faces that intersect at infinity in the case of an ideal vertex intersect in $\mathbb{H}^{3}$ if and only if their intersections with the link of the ideal vertex themselves intersect.
\end{enumerate}
Hence, by (3), we have that $\Pi_{2}$ is disjoint from the plane that truncates the vertex $A$. When $r=2$, the planes $\Pi_{1}$ and $\Pi_{2}$ will intersect if and only if their intersections with the link of $C'$ themselves intersect (by (4)). We will analyze the $r=2$ case in a moment. When $r \geq 3$, we also have that $\Pi_{2}$ does not intersect the plane that truncates the vertex $C'$, reasoned as follows. We will always choose the ``inward'' normal direction for a plane that contains a face of $T$ by indicating the appropriate opposite vertex to that face in any of our diagrams. When $r=3$, we observe that $\Pi_{2}$ contains the face of the tetrahedron (not pictured in the figure) that is the reflection of  $ABDC'$ through the face $BDC'$, and so $\Pi_{2}$ does not intersect the truncating plane of $C'$ in this case (by (3)). When $r\geq4$, then we consider the line containing the segment $BD$ which divides $\Pi_{2}$. The half of $\Pi_{2}$ that meets $C$ is prevented from intersecting the truncating plane for $C'$ by the plane $ABD$, and the other half of $\Pi_{2}$ is prevented from intersecting the truncating plane at $C'$ by the plane containing the reflection of $ABD$ through the face $BDC'$ (both of these follow from (3)). 

Therefore, when $r \neq 2$, we have that $\Pi_{2}$ has no intersection with the planes that truncate the vertices $A$ and $C'$. We observe now that these truncating planes at $A$ and $C'$ determine an open ball (i.e., the region between them in $\mathbb{H}^{3}$) which contains $\Pi_{2}$. We also note that the edge from $A$ to $C'$ is the only segment of the line of intersection of $\Pi_{1}$ with the planes $ABC'$ and $ADC'$  that lies in this ball. Using the convention for the inward normal direction given above, we conclude that, in order for $\Pi_{1}$ to intersect $\Pi_{2}$, it is necessary for that intersection to occur on the \emph{outward} side of either $ABC'$ (where inward is relative to $D$)  or the \emph{outward} side of $ADC'$ (where inward is relative to $B$), and consequently that $\Pi_{2}$ must cross at least one of the planes $ABC'$ or $ADC'$. 

By considering the geometry of the generalized vertex $B$, we have that $\Pi_{2}$ meets $ABC'$ if and only if $r=2$, and so we analyze this case now. In this case, $\Pi_{2} = BDC'$ (as planes) and $\Pi_{1}$ and $\Pi_{2}$ intersect if and only if their intersections with the link of $C'$ intersect (by (5)). The conditions for this intersection in the link of $C'$ are either $m=2$ (not possible, since $r=2$), or $n=2$ and one of $q$ or $m$ equals 3 (not possible, since $r=2$), or else $q=2$. In the last case, the intersection of $\Pi_{1}$ and $\Pi_{2}$ occurs along the edge $C'D$ at an angle of $\pi/n$, and because $q=2=r$ we must have $m \geq 6$, $p \geq 6$ and $n \geq 3$. In this case, $T=T[3,m,2;n,p,2]$ contains an immersed $(m,n,p)$ turnover, and this tetrahedron (and the set of conditions on $m, n$ and $p$) is isometric to the tetrahedron $T[2,p,n;2,m,3]$, which appears in item (1) at the end of the paper (it is listed as item (3), additionally). The summary at the end of the paper gives exact conditions on the arrangements of $l,m,q,n,p$ and $r$ which yield isometric tetrahedra. 

Otherwise, we have that $\Pi_{2}$ must intersect $ADC'$, and any possible intersection of $\Pi_{1}$ and $\Pi_{2}$ must occur on the outward side of $ADC'$ (that is, the side opposite to vertex $B$). Using the geometry of the generalized vertex $D$, we conclude that either $r=2$ (the case we just analyzed), or $p=2$, or $n=2$ and one of $p$ or $r$ equals 3. If $p=2$, then $q \geq 6$ (using the vertex $A$), $n \geq 3$ (using the vertex $D$), and $ADC'=ACDC'$ (as planes). By item (2) above, the lines $AC'$ and $CD$ are disjoint lines in the plane $ACDC'$. These lines are also the intersections with $ACDC'$ of $\Pi_{1}$ and $\Pi_{2}$, respectively. We consider the side of $ACDC'$ that is outward from vertex $B$, and the interior angles of intersection $(q-2)\pi/q$ (formed by $\Pi_{1}$ and $ACDC'$)  and $(n-1)\pi/n$ (formed by $\Pi_{2}$ and $ACDC'$) on this side of $ACDC'$ (that is, the two angles of intersection contained on this side of $ACDC'$ and in the same complementary component of these three planes).  The conditions on $n$ and $p$ imply that $(n-1)\pi/n + (q-2)\pi/q \geq \pi$, and because it is only possible for $\Pi_{1}$ and $\Pi_{2}$ to intersect to the outward side of $ACDC'$ (relative to the inward $B$ direction), we use item (4) above to  conclude that $\Pi_{1} \cap \Pi_{2} = \emptyset$ in this case. 

\begin{figure}
	\centering
	\def\svgwidth{4.5in}
	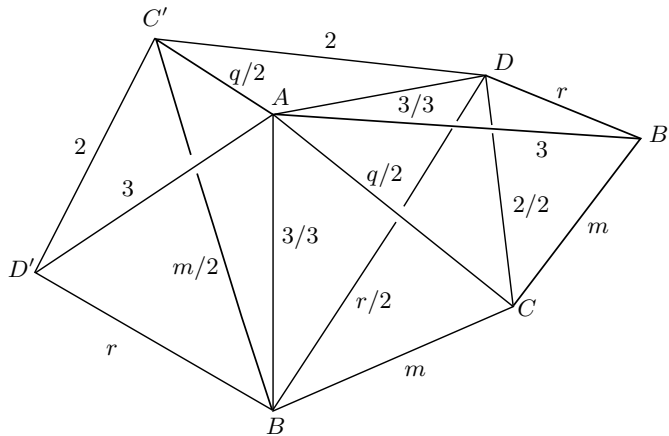
	\caption{Four copies of the tetrahedron $T[3,m,q;2,3,r]$}
	\label{F:Tet1EdgeOrder3EndGame}
\end{figure}

In the remaining case, we have $n=2$ and one of $p$ or $r$ equals 3. If $p=3$, then $r \geq 6$ and $q$ and $m$ must both be bigger than 2 and also satisfy $1/q + 1/m \leq 1/2$. We modify Figure \ref{F:Tet1EdgeOrder3} by adjoining another copy of $T$ to the face $ACD$. See Figure \ref{F:Tet1EdgeOrder3EndGame}.
In this case, $ADC'=AB'DC'$ as planes, and we consider, as in the previous case, the interior angles of intersection $(q-2)\pi/q \geq \pi/3$ and $(r-1)\pi/r \geq 5\pi/6$ formed by $AB'DC'$ with $\Pi_{1}$ and $\Pi_{2}$, respectively, on the outward side of this plane (again, relative to the inward $B$ direction). Since $(r-1)\pi/r + (q-2)\pi/q > \pi$, and again because $\Pi_{1}$ and $\Pi_{2}$ can only intersect on the side of $AB'DC'$ opposite to $B$, we conclude that $\Pi_{1} \cap \Pi_{2} = \emptyset$ in this case. The case when $n=2$ and $r=3$ is entirely similar, with the same conclusion.

\subsubsection{The single separating edge has order greater than 3:}\label{SSS:SingleEdgeCrossedOrder4OrMore} We now handle the analogous cases to the previous two cases, that is, when the search for an immersed turnover crosses a single edge between the planes $\Pi_{1}$ and $\Pi_{2}$, and when $l > 3$ (we will specify these planes in each example below, in an analogous way to the previous cases). We will show that no immersed turnovers can be found when $l > 3$. 

\begin{figure}
	\centering
	\def\svgwidth{4.5in}
	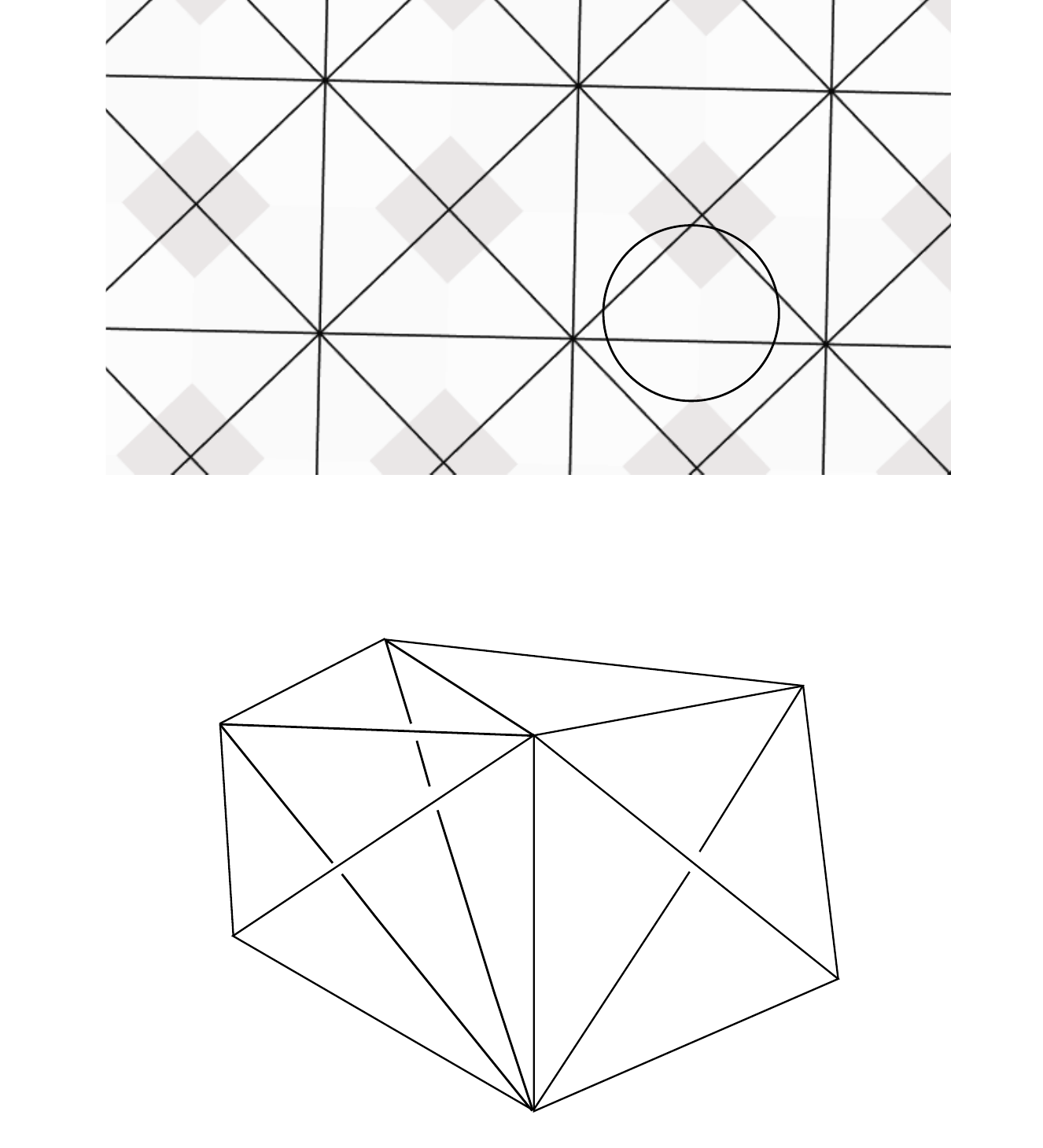
	\caption{The view from the ideal vertex of type $(2,4,4)$.}
	\label{F:RuleOut244_1}
\end{figure}

We consider first the case when $l=4$ and the vertex $B$ has the Euclidean type $(2,4,4)$ with $m=4$. See Figure \ref{F:RuleOut244_1} (we will, for the most part, drop references to the ``link'' of a vertex for the remainder of the paper, and assume that work done in, and figures referring to, the link of a vertex will be clear from the context). Referring to the lower half of this figure, we have $e_{1} = AC''$, $\Pi_{1}=AC''D'$, $e_{2}=BC$ and $\Pi_{2}=BCD$. The upper half of Figure \ref{F:RuleOut244_1} illustrates the view in the upper half-space model of $\mathbb{H}^{3}$ from the vertex $B$ (which we have placed at the point at infinity). Now $\Pi_{2}$ is represented in this diagram by the line $CD$, and the plane $\Pi_{1}$ must be represented by a circle (the circle is the boundary of a hemisphere in this model of $\mathbb{H}^{3}$). We claim that the circle representing $\Pi_{1}$ must be centered at some point in the triangle $AC''D'$, and that none of the three points $A$, $C''$ or $D'$ can be contained in this circle's interior. To see this, suppose first that the vertex $A$ of the tetrahedron is a truncated vertex. Then the plane truncating that vertex would appear as a circle in the figure. This circle would have to be centered at the point labeled $A$  because the geodesic edge from $B$ to this plane must meet the plane perpendicularly. Next, we observe that the circle representing $\Pi_{1}$ must intersect the circle centered at $A$ at a right angle (because $\Pi_{1}$ intersects the plane that truncates the vertex $A$ perpendicularly). This is only possible if the point labeled $A$ lies outside of the circle representing $\Pi_{1}$. In the case when the vertex $A$ of $T$ is an ideal vertex, then the circle representing $\Pi_{1}$ would pass \emph{through} the point labeled $A$. Since all of the vertices $A$, $C''$ and $D'$ of the tetrahedron are non-finite, the circle representing $\Pi_{1}$ cannot contain the vertices of the triangle $AC''D'$ in its interior disk. Moreover, this circle must meet each line segment $AD'$, $AC''$ and $C''D'$ (at angles of $\pi/p$, $\pi/q$ and $\pi/n$, respectively) and so the center of this circle must be contained in the triangle $AC''D'$. Such a circle is depicted. Since any such circle cannot intersect the line $CD$, we conclude that $\Pi_{1} \cap \Pi_{2} = \emptyset$. An analogous argument can be used to show that we obtain no immersed turnover in this fashion, whenever the vertex $B$ is Euclidean and $l$ is not equal to  2 or 3 (this occurs only when the triple $(l,m,r)$ equals one of $(4,2,4)$, $(6,2,3)$ or $(6,3,2)$).

\begin{figure}
	\centering
	\def\svgwidth{4.5in}
	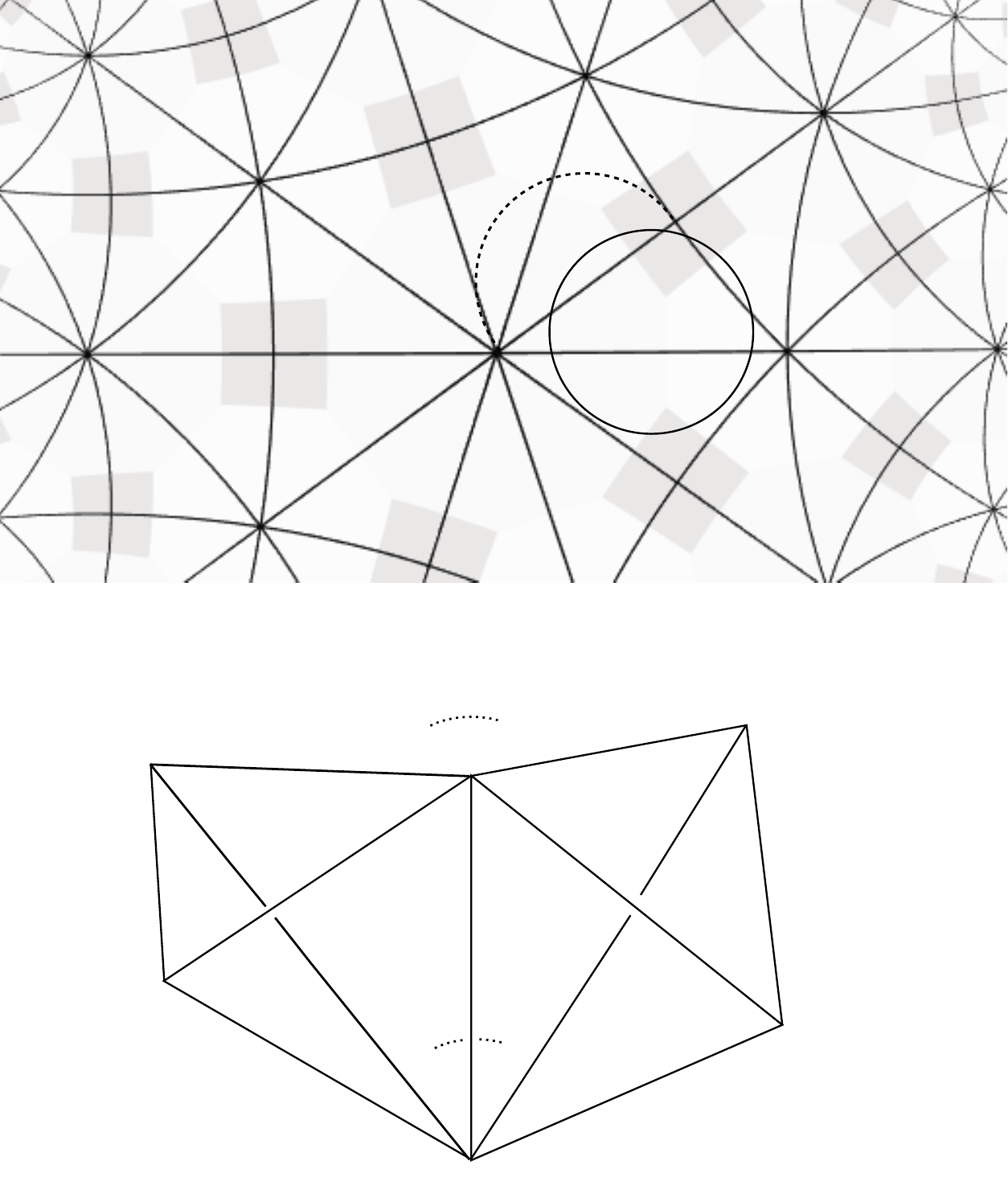
	\caption{The view from the truncated vertex of type $(2,4,5)$.}
	\label{F:RuleOut2xy_1}
\end{figure}

\label{arg:ReferToThisArgLater2}We are left then to consider the case when $l \geq 4$ and the vertex $B$ has a hyperbolic type. The argument is similar to the Euclidean vertex case, but we provide the details. Consider first the case of Figure \ref{F:RuleOut2xy_1} (the upper part of this figure, along with the similar figures in this section, was generated using the software \emph{KaleidoTile} by Jeffrey Weeks \cite{Weeks_KaleidoTile}). For the purposes of illustration, we have assumed that the vertex $B$ has the type $(2,4,5)$, with $l=5$, $m=2$ and $r=4$. Here, we have $e_{1}=AD''$, $\Pi_{1}=AD''C''$, $e_{2}=BC$ and $\Pi_{2}=BCD$. We consider the hyperbolic plane $\Pi_{B}$ that truncates vertex $B$ as a hemisphere in the upper half-plane model, and wish to construct a ``view from $B$'' that is similar to the previous case when the $B$ was an ideal vertex. The Poincar\'e disk $(2,4,5)$ tiling pattern of the figure results from projecting this \emph{hemi}sphere to the bounding plane of $\mathbb{H}^{3}$ through the south pole of the \emph{whole} sphere that contains it \cite[Figure 2.12, p. 58]{Thurston97}. An important observation about this projection is that it is equivalent to projecting every point $x \in \Pi_{B}$ to the bounding plane of half-space along the geodesic ray that is perpendicular to $\Pi_{B}$ at $x$.  In particular, as in the Euclidean vertex case, each line or circular arc in the figure is the ideal boundary of a plane (each plane corresponding to a face in the tiling of $\mathbb{H}^{3}$ by $T$) that meets  $\Pi_{B}$ perpendicularly, and this projection is conformal, so that the angle of intersection between two lines or circular arcs in the figure is equal to the angle of intersection of the corresponding planes in $\mathbb{H}^{3}$. We have indicated, in the projection of the figure, the images of the intersection of five copies of $T$ with $\Pi_{B}$, labeled the endpoints of the lines emanating from $B$ by the corresponding letters in the lower part of the figure, and applied an isometry so that $A$ (or, in the case that the vertex $A$ is truncated, the center of the circle that represents the truncating plane for the vertex $A$) is at the center of the Poincar\'e disk. The planes $\Pi_{1}$ and $\Pi_{2}$ are represented by a circle and the circular arc $CD$, respectively.  
We observe that, if the vertex $C''$ is truncated, then the truncating plane $\Pi_{C''}$ for $C''$ will appear in the figure as a circle (not pictured) with center on the segment $AC''$, because the point labeled $C''$ is the endpoint of a semicircle in the half-space model that is perpendicular to both $\Pi_{B}$ and $\Pi_{C''}$ (to see this, recall that we may consider the projection from $\Pi_{B}$ to the bounding plane as a projection along arcs of such semicircles). As in the previous case, the point $C''$ cannot be contained in the interior of the circle that is the ideal boundary of $\Pi_{1}$, because then the arc of the semicircle from $C''$ to its inverse image in $\Pi_{B}$ under the projection would meet $\Pi_{1}$, and this is impossible because this arc meets $\Pi_{C''}$ perpendicularly and $\Pi_{C''}$ and $\Pi_{1}$ are orthogonal (the contradiction arises because it would imply the existence of a triangle with two right angles). The same argument holds when either of $A$ or $D''$ is a truncated vertex, and therefore, as in the previous case, we have that the ideal boundary of $\Pi_{1}$ must bound a disk whose interior is disjoint from the points $A$, $C''$ and $D''$ (these points may lie on the ideal boundary of $\Pi_{1}$ if they are ideal vertices of $T$). The ideal boundary of $\Pi_{1}$ intersects the segments $AC''$ and $AD''$ and the circular arc $C''D''$ at angles of $\pi/q$, $\pi/p$ and $\pi/n$, respectively, and the center of the circle representing this ideal boundary has its center contained in the hyperbolic triangle $AC''D''$ in the projection. This is the circle that is depicted in the figure. But such a circle can have no points in the hyperbolic polygon $CAD''C''D'C'D$ that lie outside of the union of hyperbolic  triangle $AC''D''$ and the circle with the segment $AC''$ as its diameter (pictured with a dashed arc in the figure). Consequently, this circle cannot meet any of the sides of this hyperbolic polygon other than $AD''$ and $D''C''$, and, in particular, we have $\Pi_{1} \cap \Pi_{2} = \emptyset$. An analogous argument works whenever $B$ has hyperbolic type with one incident order 2 edge and $\l \geq 4$.

\begin{figure}
	\centering
	\def\svgwidth{4.5in}
	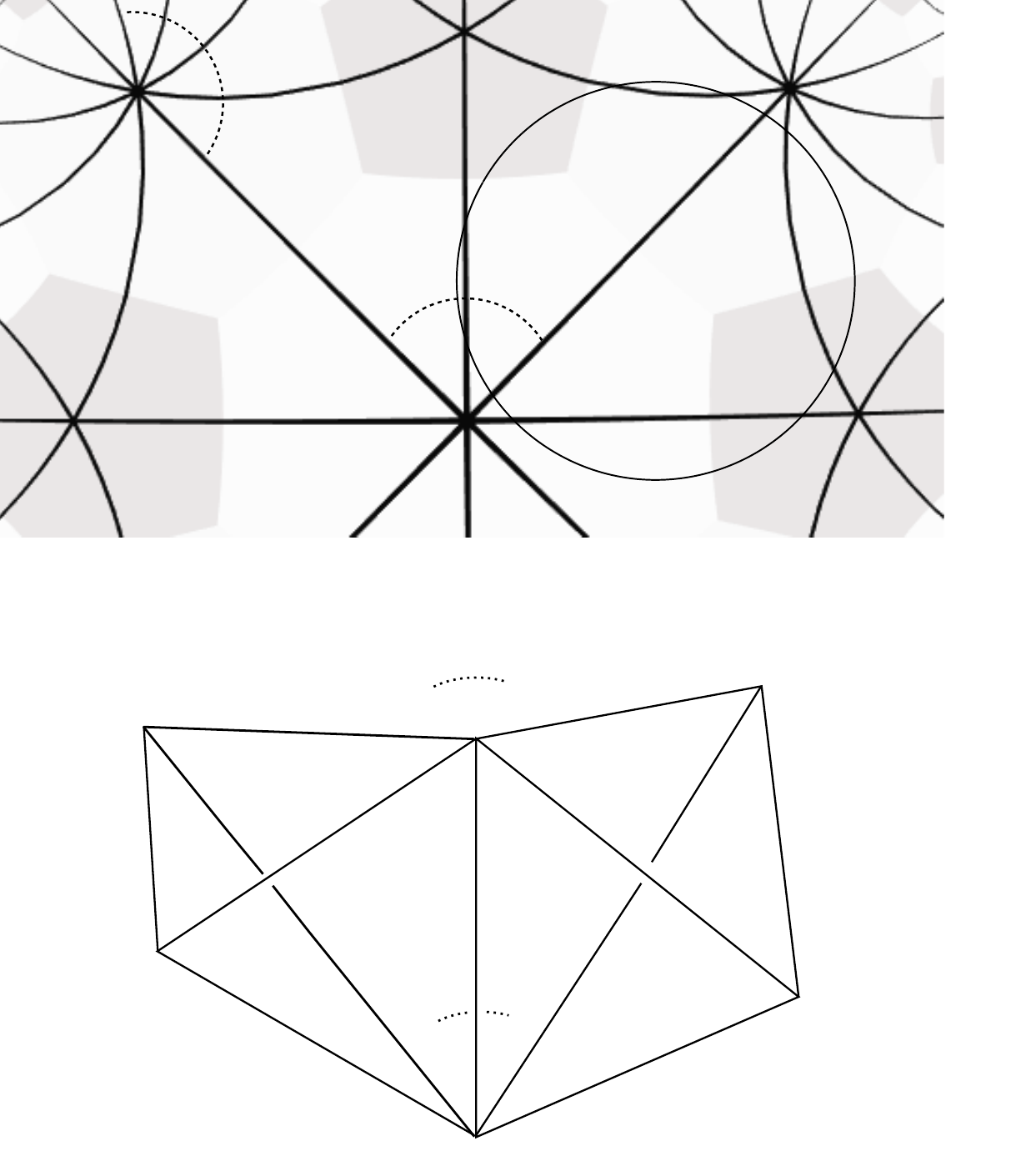
	\caption{The view from the truncated vertex of generic hyperbolic type when none of $l$, $m$ and $r$ is 2.}
	\label{F:RuleOutxyz_1}
\end{figure}

\begin{figure}
	\centering
	\def\svgwidth{4.5in}
	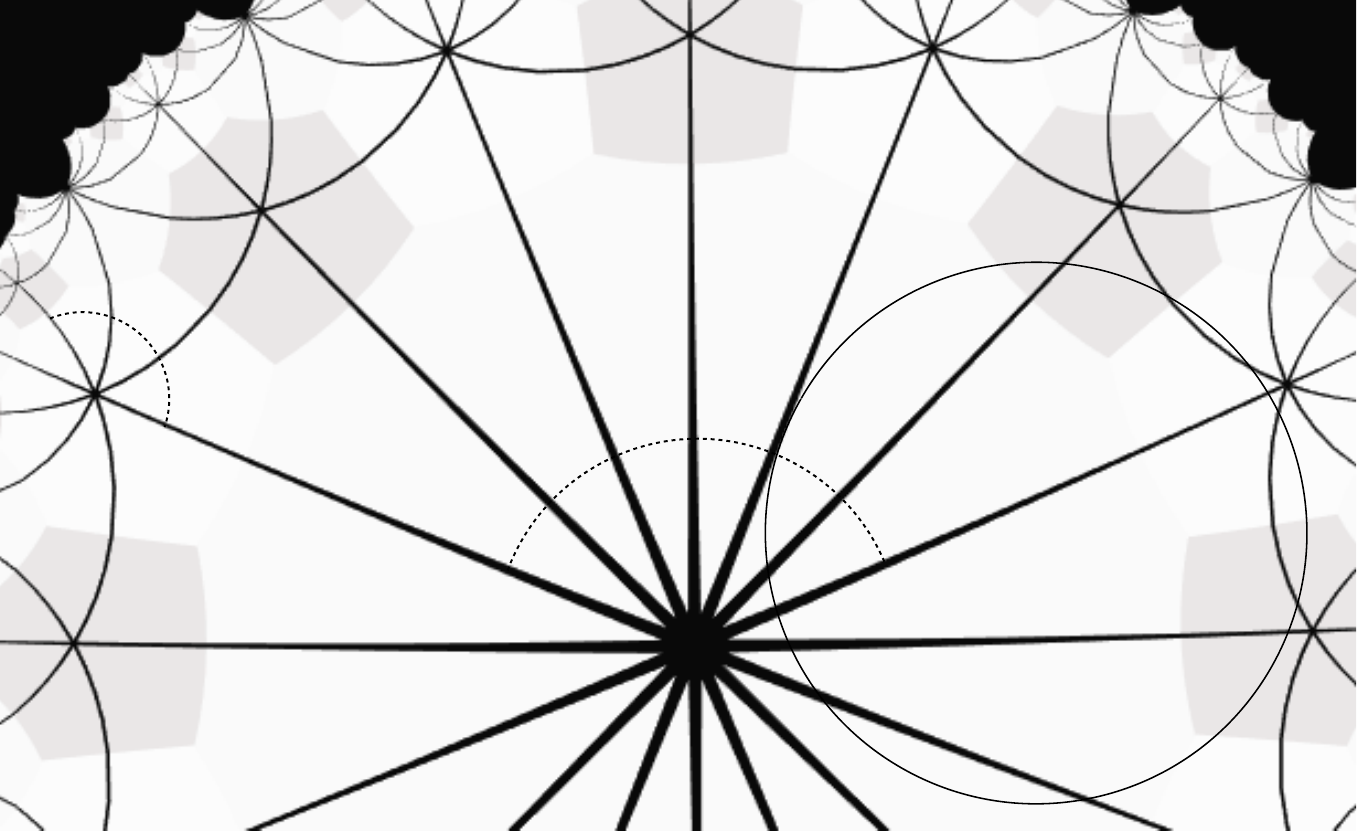
	\caption{Another view from the truncated vertex of generic hyperbolic type when none of $l$, $m$ and $r$ is 2.}
	\label{F:RuleOutxyz_2}
\end{figure}

The case when $\l \geq 4$ and $B$ has hyperbolic type with no incident order 2 edge is similar. See Figure \ref{F:RuleOutxyz_1}, in which $\Pi_{2}$ is represented by the circular arc $CD$ and $\Pi_{1}$ is represented as the circle pictured. 
When $l \geq 4$, we observe that, in any similar picture (for example, Figure \ref{F:RuleOutxyz_2}), the angles $\alpha=(l-2)\pi/l$ and $\beta=(r-1)\pi/r$ will always be at least $\pi/2$.
Consequently, since $\alpha \geq \pi/2$ and because the center of the circle representing $\Pi_{1}$ is contained in the hyperbolic triangle $AEF$, this circle will be disjoint from the interior of the segment $AD$ (it may pass through $A$, if the corresponding vertex is ideal). Also, noting that $AD$ will always have Euclidean length equal to one of the lengths $|AF|$ or $|AE|$, the conditions on $\alpha$ and $\beta$ imply that no point of the circle $CD$ that lies above the line $AD$ will be closer to the center of the circle representing $\Pi_{1}$ than any of the points $A$, $E$ or $F$. Since $A$, $E$ and $F$ are not contained in the interior of this circle, we can conclude that $\Pi_{1} \cap \Pi_{2} = \emptyset$ in this case.

\subsection{The case when multiple edges separates $e_{1}$ from $e_{2}$:}\label{SS:MultipleEdgesCrossed} 

Recall that $\Pi_{\cal{T}}$ denotes the plane stabilized by a copy of a triangle subgroup in the fundamental group of $\cal{O}_{T}$, and that $e_{1}$ and $e_{2}$ denote two developed coplanar edges of $T$ whose (perpendicular) intersections with $\Pi_{\cal{T}}$ correspond to two of the cone points of an immersed turnover (whose fundamental group is the triangle group stabilizing $\Pi_{\cal{T}}$) in $\cal{O}_{T}$. 

\begin{figure}
	\centering
	\def\svgwidth{6.5in}
	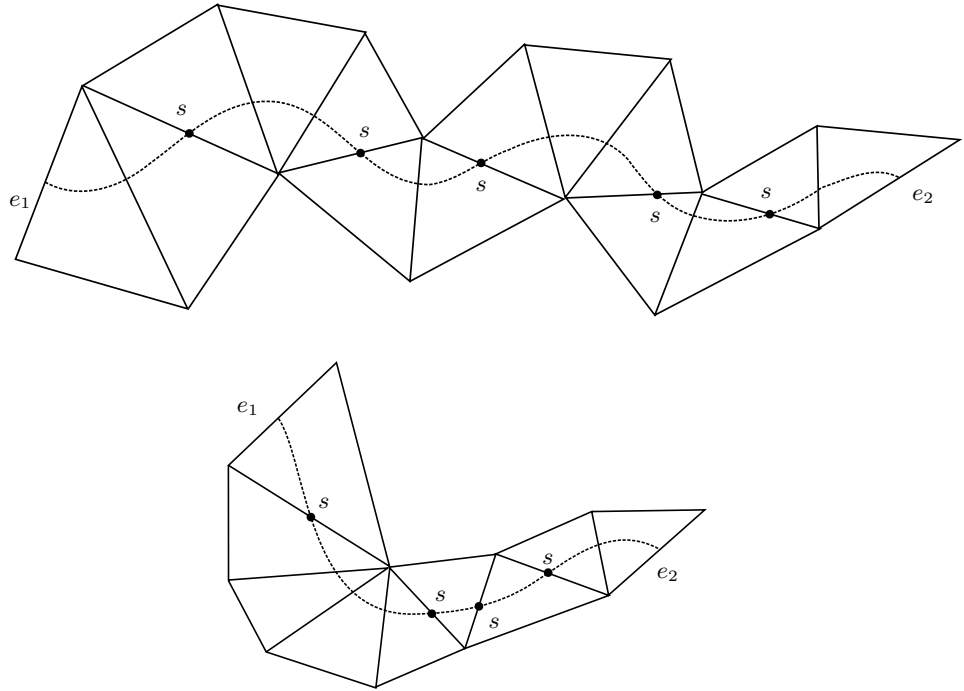
	\caption{Schematic of some possible developments of a face of $T$, together with switches and the possible intersection of the plane $\Pi_{\cal{T}}$.}
	\label{F:FaceDevelopment}
\end{figure}

\begin{figure}
	\centering
	\def\svgwidth{4in}
	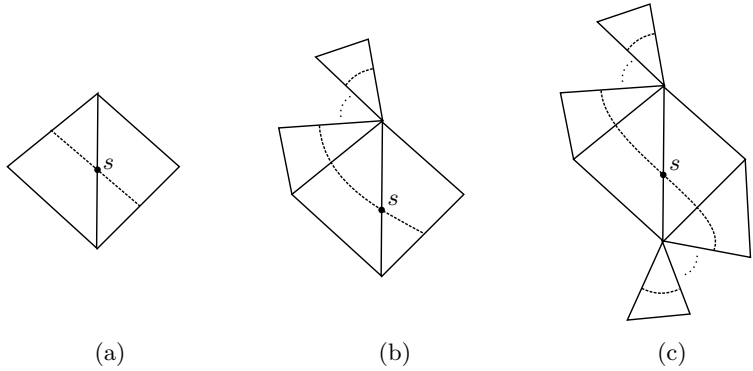
	\caption{One type of possibility for the subdiagram components for a diagram of the type given in Figure \ref{F:FaceDevelopment}.}
	\label{F:FaceSubdiagrams1}
\end{figure}

\begin{figure}
	\centering
	\def\svgwidth{4in}
	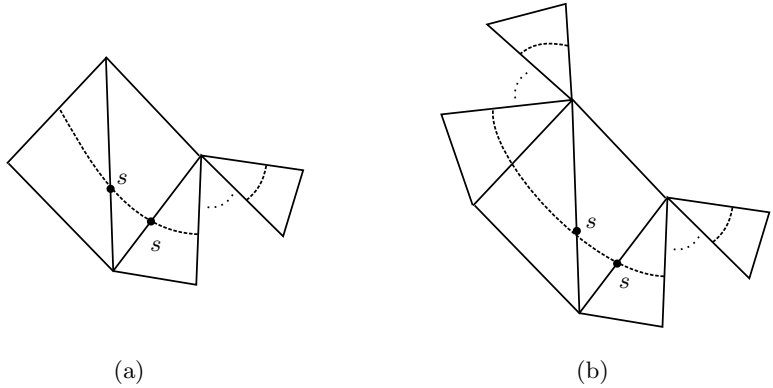
	\caption{Another type of possibility for the subdiagram components for a diagram of the type given in Figure \ref{F:FaceDevelopment}.}
	\label{F:FaceSubdiagrams2}
\end{figure}

\emph{Notation.} For the remainder of the paper, $\Pi_{F}$ refers to the plane containing $e_{1}$ and $e_{2}$. It is the development in $\mathbb{H}^{3}$ of one face $F$ of $T$. The diagrams from Figures \ref{F:FaceDevelopment}, \ref{F:FaceSubdiagrams1} and \ref{F:FaceSubdiagrams2} (along with several other figures later in this section) are all drawn with the convention that $\Pi_{F}$ is the page containing the illustration. We use $L_{F}$ to denote the intersection of $\Pi_{\cal{T}}$ with $\Pi_{F}$. Additionally, the phrase ``the other side of $\Pi_{F}$'' refers, in each of the relevant figures, to the side of $\Pi_{F}$ that is behind the page (relative to the reader), and the use of the word ``plane'' at any edge in a diagram \emph{always} refers to a plane that is the development of a face of $T$ in $\mathbb{H}^{3}$ that passes through that edge.

Now that we have determined the conditions on $T$ which give rise to a turnover in $\cal{O}_{T}$ when $\Pi_{\cal{T}}$ intersects a single edge in the development of $F$ between $e_{1}$ and $e_{2}$, we will show that it is impossible for there to be more than one such edge in the development of $F$ between $e_{1}$ and $e_{2}$. This will complete the classification of immersed turnovers in tetrahedral orbifolds with no finite generalized vertices.

Figure \ref{F:FaceDevelopment} shows two possible schematic diagrams for this discussion. In each of the subfigures, the edges $e_{1}$ and $e_{2}$ are indicated, and the dotted line represents $L_{F}$. Notice that, in each triangle of the planar development of $F$, there is always a unique translate of a vertex of $T$ that is separated from the other two by $\Pi_{\cal{T}}$. The edge translates of $T$ labeled by ``$s$'' represent points at which this vertex switches.
We consider the following procedure for dividing any diagram of the type from Figure \ref{F:FaceDevelopment} into subdiagrams of the type (up to possible reflection or order two rotation) given in Figures \ref{F:FaceSubdiagrams1} and \ref{F:FaceSubdiagrams2}: 
\begin{enumerate}
	\item Starting at the first edge of the diagram, we follow $L_{F}$ until we arrive at the first switch. There must always be such a switch, for otherwise the supposed turnover would be parallel to a cover of an embedded turnover corresponding to one of the truncated vertices of $T$.
	\item If the switch is the only switch in the diagram, then our diagram looks like, up to reflection or rotation, one of the diagrams from Figure \ref{F:FaceSubdiagrams1}. In this case, we stop.
	\item If there is more than one switch and the diagram looks like, up to reflection or rotation, one of the diagrams from Figure \ref{F:FaceSubdiagrams2}, then we stop.
	\item If we have not halted in the previous two steps, then the diagram up to and including the first edge after the first switch looks like the diagram in either Figure \ref{F:FaceSubdiagrams1}(a) or \ref{F:FaceSubdiagrams1}(b). Call this portion a \emph{subdiagram}. 
	\item Starting at the last edge of the subdiagram from the previous step, we repeat this process with the remaining portion of the original diagram, starting from the first step, until we reach edge $e_{2}$.
\end{enumerate}
This procedure divides our diagram into subdiagrams of the type illustrated in Figure \ref{F:FaceSubdiagrams1}(a) and \ref{F:FaceSubdiagrams1}(b), with the possible exception that the final subdiagram may be of the type in Figure \ref{F:FaceSubdiagrams1}(c) or one of the two types in Figure \ref{F:FaceSubdiagrams2} (we note that this process can eliminate certain switches, in each of the resulting subdiagrams). Again, we denote by $\Pi_{1}$ and $\Pi_{2}$ the planes at $e_{1}$ and $e_{2}$, respectively, whose intersections with $\Pi_{\cal{T}}$ are supposed to form two of the sides of a triangle in the tiling of $\Pi_{\cal{T}}$. Our strategy is to use the subdiagrams of Figures \ref{F:FaceSubdiagrams1} and \ref{F:FaceSubdiagrams2} to find a sequence of planes in $\mathbb{H}^{3}$---one or more planes at each of the two outer-most edges of each subdiagram---that are pairwise disjoint on either side of $\Pi_{F}$ and that therefore separate $\Pi_{1}$ from $\Pi_{2}$.  

\begin{figure}
	\centering
	\def\svgwidth{5in}
	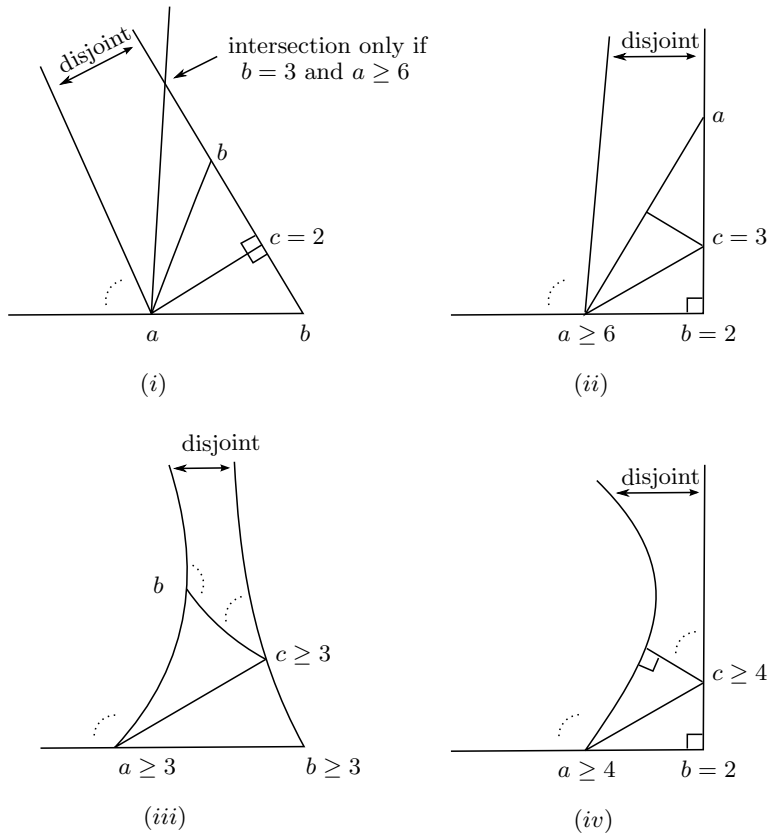
	\caption{The possibilities for the intersection of lines in a triangular tiling of $\mathbb{E}^{2}$ or $\mathbb{H}^{2}$.}
	\label{F:TriangleLineIntersections}
\end{figure}

We first make two observations about the subdiagram from Figure \ref{F:FaceSubdiagrams1}(a). First, if either of the orders of the two edges separated by the switch  is 2, then no plane at either edge can meet any of the planes at the other edge (excepting the plane $\Pi_{F}$). This fact follows from the extensive analysis done in Subsection \ref{SS:SingleEdgeCrossed}. Second, if the two planes at the outer edges that are inclined closest toward the switch (``inclined closest''  means closest, on the other side of $\Pi_{F}$, to the planes that pass through the switch edge) do meet (thus generating an immersed turnover in $\cal{O}_{T}$ with two singular points of order at least 6 and one singular point of order at least 3), then the next two planes (one at either outer edge) inclined away from the switch do not meet. This fact follows from an easy analysis of the patterns of line intersections in hyperbolic triangular tilings. See Figure \ref{F:TriangleLineIntersections} for the conditions on the vertex orders of an $(a,b,c)$ hyperbolic triangular tiling under which such intersections can occur. In this case (although this will not be the case for subsequent applications of this figure), Figure \ref{F:TriangleLineIntersections} should be thought of as depicting the plane which meets $\Pi_{F}$ and the two southwest-to-northeast edges from Figure \ref{F:FaceSubdiagrams1}(a) perpendicularly, so that all the planes through these two edges appear as lines in Figure \ref{F:TriangleLineIntersections}.  In particular, in order for the next two planes inclined away from the switch in Figure \ref{F:FaceSubdiagrams1}(a) to meet, then one of the southwest-to-northeast edges must have order 2, which does not happen in this situation.

Therefore, it is left to show that, for each of the remaining types of subdiagram, the two planes at the outer-most edges that are inclined closest to the single or double switch in the subdiagram do not intersect (again, ``inclined closest'' means closest, on the other side of $\Pi_{F}$, to the planes passing through the switch edge(s)). This will produce the sequence of planes that separates $\Pi_{1}$ and $\Pi_{2}$, and therefore complete the proof. We will show this by cases, which are indicated by their labels in the figures.

\subsubsection{\ref{F:FaceSubdiagrams1}(b):}\label{SSS:FaceSubdiagrams1b} 
\begin{figure}
	\centering
	\def\svgwidth{3in}
	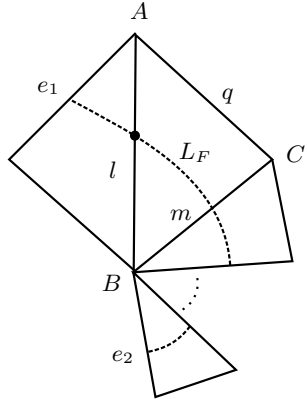
	\caption{The case of Figure \ref{F:FaceSubdiagrams1}(b).}
	\label{F:OneSwitch1}
\end{figure}

See Figure \ref{F:OneSwitch1}, in which we have supposed without loss of generality that $F$ is the face $ABC$ of the tetrahedron $T$, as in Figure \ref{F:TLMQNPR}. This picture only differs from Figure \ref{F:FaceSubdiagrams1}(b) by a 180$^\circ$ rotation. Observe that the edges incident at the vertices $A$ and $B$ have orders $l$, $q$, $p$ and $l$, $m$, $r$ (respectively). 

\begin{figure}
	\centering
	\def\svgwidth{5in}
	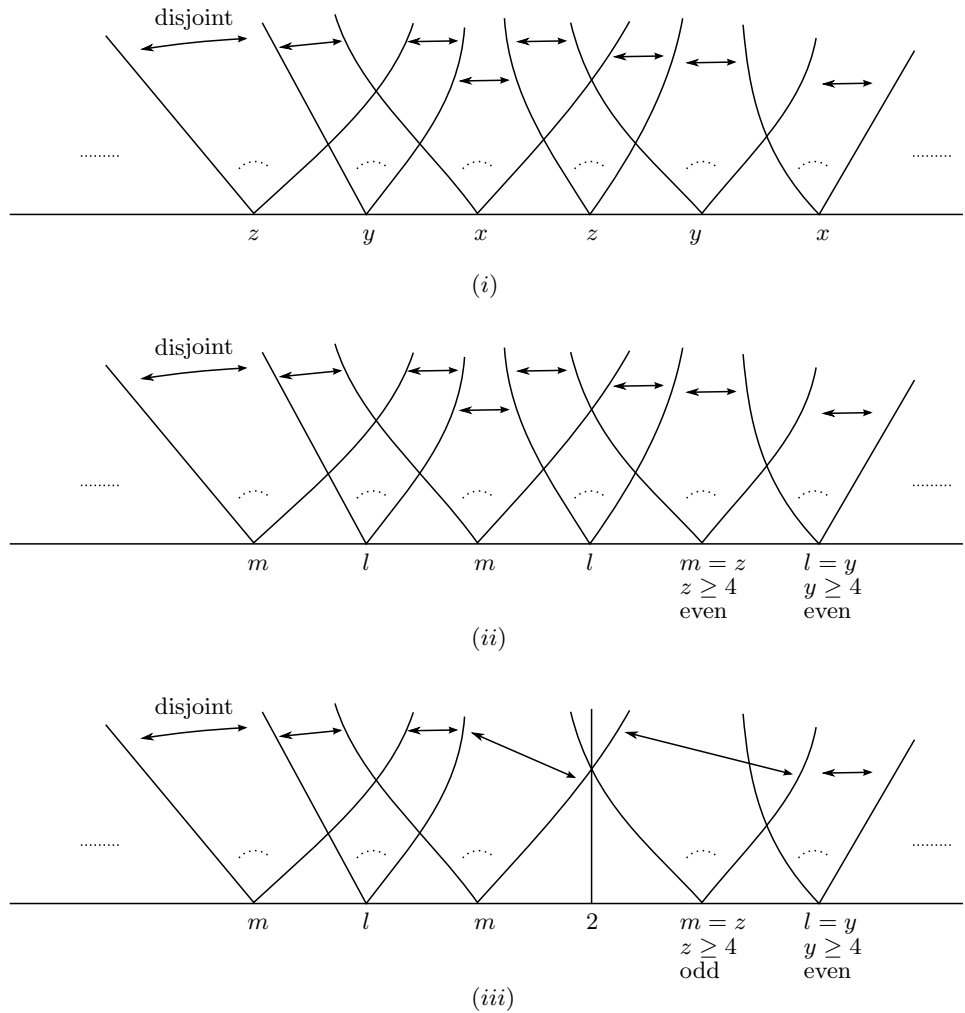
	\caption{Patterns of intersections of certain lines corresponding to sides in a triangular tiling of $\mathbb{H}^{2}$ or $\mathbb{E}^{2}$. Double arrows indicate two lines that do not intersect above the horizontal line.}
	\label{F:UseLineIntersectionsCase1}
\end{figure}

We observe that the vertex $B$ must have at least one order 2 edge incident to it. Otherwise, if $B$ were of the type $(x,y,z)$ with all orders at least 3, then it is readily seen, by using the information from Figure \ref{F:TriangleLineIntersections}$(iii)$ applied to vertex $B$, that $\Pi_{2}$ (the plane through $e_{2}$ that is inclined closest to the switch) cannot meet the plane at edge $BC$ that is inclined closest to the switch. We indicate how this can be determined. Recall that we may construct the view from $B$ as a triangular tiling of either the Euclidean or hyperbolic plane (in this case, a tiling by $(x,y,z)$ triangles) such that $\Pi_{F}$ appears as a horizontal line, and such that each edge incident to $B$ appears as a point on that line and each plane through an edge incident to $B$ appears as a line (or hyperbolic line, if $B$ is super-ideal) passing through the corresponding point in the view from $B$. Using Figure \ref{F:TriangleLineIntersections}$(iii)$, we can conclude that the view from $B$, when $B$ has no incident order 2 edge, looks schematically like Figure \ref{F:UseLineIntersectionsCase1}$(i)$. 
This figure assumes that $x$, $y$ and $z$ are all odd; the other cases are similar. Suppose, for example, that the right-most point $x$ in this figure represents the edge $BC$ ($x$ also indicates the order of that edge), and that the (schematic) line through this point inclined furthest to the right represents the plane through edge $BC$ inclined closest to the switch edge $AB$. Then it is easily seen that no right-most inclined line through any subsequent point to the left  along the horizontal can intersect with this line. Consequently, the planes to which these lines correspond cannot intersect on the other side of $\Pi_{F}$ (i.e., the other side of the page in Figure \ref{F:OneSwitch1}). In particular, $\Pi_{2}$ cannot cross the plane through $BC$ inclined closest to the switch, as we wished to show. Furthermore, by our analysis in the cases of Subsection \ref{SS:SingleEdgeCrossed}, the only way that $\Pi_{1}$ can meet the plane through edge $BC$ that is inclined closest to the switch is if $B$ has an incident order 2 edge. Consequently, if there is no such order 2 edge at $B$, then we have $\Pi_{1} \cap \Pi_{2} = \emptyset$. 

\begin{figure}
	\centering
	\def\svgwidth{5in}
	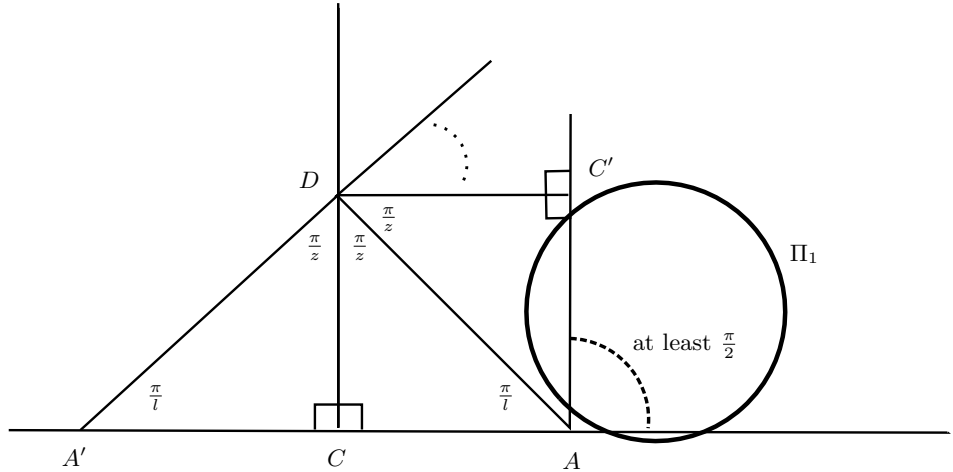
	\caption{The schematic view from the vertex $B$ in the case when $l =y \geq 4$, $m=2$ and $r=z \geq 4$.}
	\label{F:OneSwitch1RuleOut1}
\end{figure}

So $B$ either has the type $(2,3,x\geq6)$ or $(2,y\geq4,z\geq4)$. In the latter case, if $l=y$ or $l=z$, then we have shown in Subsection \ref{SSS:SingleEdgeCrossedOrder4OrMore} that $\Pi_{1}$ is disjoint from every plane through edge $BC$. If $l=y$ and $m=z$ and $l$ and $m$ are both even, then it is a simple exercise, using Figure \ref{F:TriangleLineIntersections}$(i)$, to show that no plane that is inclined closest to the switch edge $AB$ through any of the subsequent edges from $BC$ toward $e_{2}$ along $L_{F}$ can meet the plane through edge $BC$ that is inclined closest to the switch, as in the argument of the previous paragraph (the schematic of the view from $B$ in this case would be Figure \ref{F:UseLineIntersectionsCase1}$(ii)$, with the edges $AB$ and $BC$ corresponding to the right-most points labeled $l$ and $m$, respectively). So $\Pi_{1} \cap \Pi_{2} = \emptyset$ in this case. If $l=y$ and $m=z$ and $m$ is odd, we can use the same argument (this time using the information from items $(i)$, $(ii)$ and $(iv)$ from Figure \ref{F:TriangleLineIntersections} to obtain the schematic view from $B$ as depicted in Figure \ref{F:UseLineIntersectionsCase1}$(iii)$) to conclude that $\Pi_{1} \cap \Pi_{2} = \emptyset$. The analogous cases, where $l=y$ and $m=z$ and $l$ and $m$ are of mixed parity, are similar.  The case when $l=y$ or $l=z$ and $m=2$ requires a bit more analysis. \label{arg:ReferToThisArgLater1} In this case, we use the geometry of the vertex $A$, the fact that $l \geq 4$ and the information from Figure \ref{F:TriangleLineIntersections} to conclude that $\Pi_{1}$ cannot intersect the plane through edge $AB$ that is inclined closest to the edge $BC$. But $\Pi_{1}$ must intersect $\Pi_{F}$ and it must intersect some of the planes through the switch edge $AB$. We refer to Figure \ref{F:OneSwitch1RuleOut1}, which depicts the schematic view from $B$ in this case, with $l=y \geq 4$ and $m=2$ (the third edge incident to $B$, which would have the label $r$ in the tetrahedron $T$, is labeled by $z \geq 4$).  
In this figure, the line segment $AD$ corresponds to the plane through the switch edge $AB$ of the tetrahedron that is inclined closest to the edge $BC$. As we have seen in previous cases, the ideal boundary of $\Pi_{1}$, in this view, is a circle that cannot contain any vertex of the triangulation in its interior disk. Since $z \geq 4$, we may conclude from the figure that the ideal boundary of $\Pi_{1}$ cannot intersect the line $A'D$. By noting that the line $A'D$ represents the plane inclined closest to the switch through the edge just after the edge $BC$ along $L_{F}$ toward $e_{2}$ in Figure \ref{F:OneSwitch1}, we may use the previous arguments from this paragraph to conclude that $\Pi_{1} \cap \Pi_{2} = \emptyset$ in this case.

Referring to the first sentence of the previous paragraph, in the latter case and when $l=2$ and $y=4=z$, we may show that $\Pi_{1} \cap \Pi_{2} = \emptyset$ by using the Euclidean vertex argument as in Figure \ref{F:RuleOut244_1}. In the latter case and when $l=2$ and one of $y$ or $z$ is greater than 4, it is again readily shown that the second closest plane to the switch through edge $BC$ (recall that $\Pi_{1}$ must be disjoint from this plane, by the observation of the penultimate paragraph before the start of this subsection) misses the plane inclined closest to the switch at every subsequent edge that $L_{F}$ crosses toward $e_{2}$. The argument uses the information of items $(i)$, $(ii)$ and $(iv)$ from Figure \ref{F:TriangleLineIntersections}, and is similar to the arguments already presented in the previous two paragraphs. Thus, we have $\Pi_{1} \cap \Pi_{2} = \emptyset$ in the case that the type of vertex $B$ is $(2, y\geq 4, z \geq 4)$.

\begin{figure}
	\centering
	\def\svgwidth{4in}
	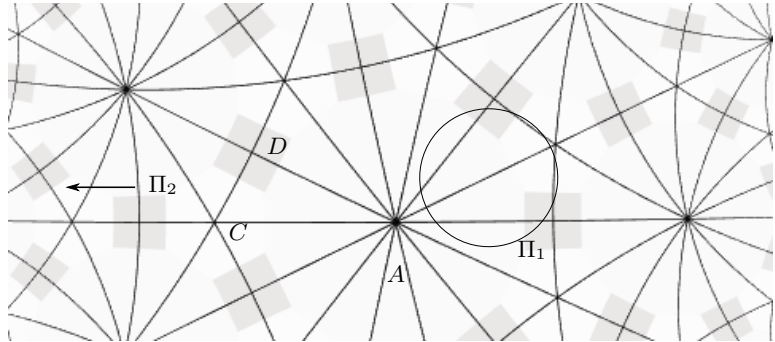
	\caption{A view from the truncated vertex of hyperbolic type $(2,3,7)$. The arrow indicates that the plane $\Pi_{2}$ is represented by a circular arc that meets the horizontal somewhere to the left of the arc $CD$.}
	\label{F:RuleOut23x_1}
\end{figure}

\label{arg:ReferToThisArgLater3}This leaves us with the possibility that $B$ has type $(2,3,x\geq6)$. When $l=x$, then we are in a case that is similar to the first case in Section \ref{SSS:SingleEdgeCrossedOrder4OrMore}, i.e., we have to consider a regular $l$--gon in either the Euclidean or hyperbolic plane and a circle centered inside the polygon that does not contain in its interior the center of the polygon, any vertex of the polygon or any midpoint of a side. In this case, however, we observed that such a circle (representing $\Pi_{1}$) must be disjoint from all but two sides of the polygon. But the plane $\Pi_{2}$ will correspond in such a picture to a line or circular arc in the picture that does not meet the interior of this polygon, and so $\Pi_{1} \cap \Pi_{2} = \emptyset$ when $l=x$. See Figure \ref{F:RuleOut23x_1} for an example illustration of this argument, in the case when $x=7$.

\begin{figure}
	\centering
	\def\svgwidth{4in}
	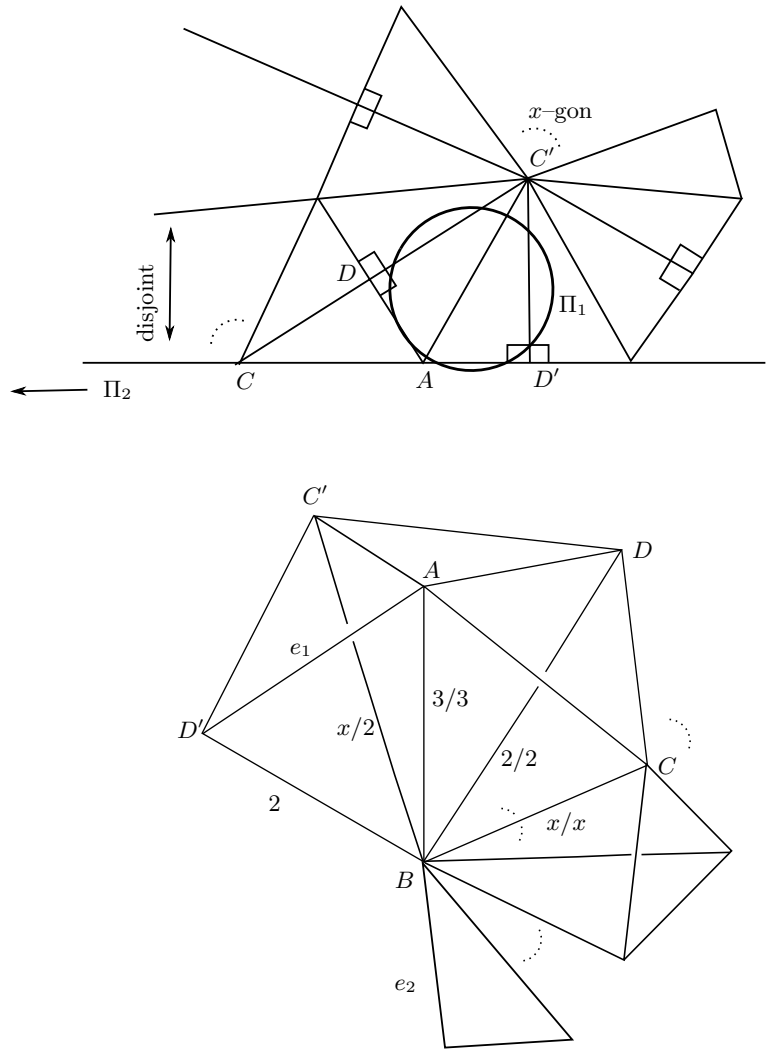
	\caption{A view from the truncated vertex of type $(2,3,x\geq6)$.}
	\label{F:RuleOut23x_2}
\end{figure}

The cases when $l=2$ or $l=3$ remain. In the case when $l=3$, we refer to Figure \ref{F:RuleOut23x_2}. 
The upper half of this figure depicts the salient aspects of the view from vertex $B$, as in the previous cases we have considered. The lower half of the figure depicts part of the development of $T$ in $\mathbb{H}^{3}$. In particular, in the lower half of the figure, the triangle with edge $e_{2}$ and the lowest set of elliptical dots are both meant to lie in $\Pi_{F}$ (which is the horizontal line $CAD'$ in the upper half of the figure), and the plane $\Pi_{2}$ is not depicted, although $\Pi_{1}=AC'D'$ is. In the upper half of the figure, $\Pi_{1}$ is represented by a circle centered at some point inside the triangle $AC'D'$ that cannot meet any vertex of the triangulation and that can only meet the sides $AD$ and $AD'$ of the $x$--gon centered at $C'$ (the fact that this circle can meet no other sides of the $x$--gon centered at $C'$ follows by an argument similar to that depicted in Figure \ref{F:RuleOut2xy_1} from Section \ref{SSS:SingleEdgeCrossedOrder4OrMore}). Since $\Pi_{2}$ must be represented by a line emanating from a vertex on the line $CAD'$ which is further to the left than $C$ (the direction, in the upper part of the figure, to which the line representing $\Pi_{2}$ must lie is indicated by the lower left arrow), and no such lines will enter the $x$--gon centered at $C'$, we conclude that $\Pi_{1} \cap \Pi_{2} = \emptyset$ in this case.

\begin{figure}
	\centering
	\def\svgwidth{4in}
	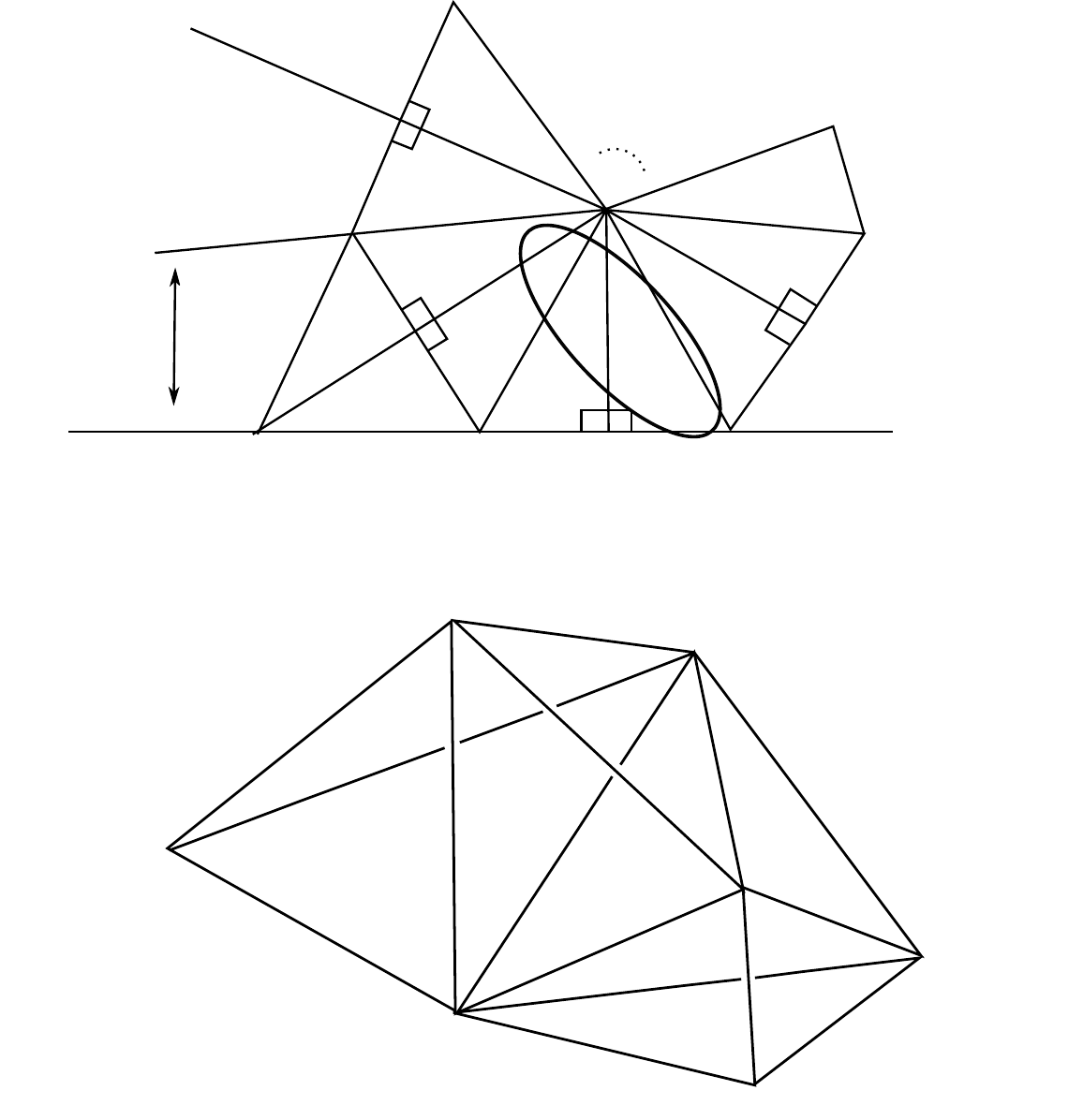
	\caption{A view from the truncated vertex of type $(2,3,x\geq6)$.}
	\label{F:RuleOut23x_3}
\end{figure}

When $l=2$, then the only way for which we are unable to apply the preceding argument is when $m=3$. See Figure \ref{F:RuleOut23x_3}. 
This is because the angle $\angle A'DC'$ is less than $\pi/2$ when $x > 6$, and so it is, in principle, possible that the circle representing $\Pi_{1}$ (whose center must be contained in the triangle $AC'D$) may intersect the line representing $\Pi_{2}$ if $\Pi_{1}=AC'D$ and $\Pi_{2}=A'D'BD$ (we have drawn the circle as an ellipse in the view from $B$ in order to indicate this possible intersection). However, using the accompanying tetrahedral illustration and the techniques of Subsection \ref{SSS:SingleEdgeCrossedOrder2} (applied to vertex $D$), it is readily seen that we must have $p=2$ and $n=3$ in order for $\Pi_{1}$ and $\Pi_{2}$ to intersect. However, because we assume that $T$ has no finite vertices and because $l=2$, we do not allow $p=2$. (\emph{Note:} When $l=p=2$ and $m=n=3$ (so that the vertex $A$ is finite), there is an immersed turnover of type $(q,x,x)$ in $T$, provided that $q\geq3$ and $x\geq4$. See the conjectural classification at the end of this paper. In this case, $T=T[2,3,q;3,2,x]$, which is isometric to the tetrahedron listed in item (6).)

\subsubsection{\ref{F:FaceSubdiagrams1}(c):}\label{SSS:FaceSubdiagrams1c} 
\begin{figure}
	\centering
	\def\svgwidth{2in}
	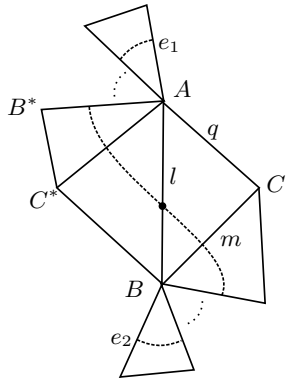
	\caption{The case of Figure \ref{F:FaceSubdiagrams1}(c).}
	\label{F:OneSwitch2}
\end{figure}

See Figure \ref{F:OneSwitch2}, in which again we have supposed without loss of generality that $F$ is the face $ABC$ of the tetrahedron $T$, with the edges incident at the vertices $A$ and $B$ having orders $l$, $q$, $p$ and $l$, $m$, $r$ (respectively).  We again denote by $\Pi_{1}$ and $\Pi_{2}$ the planes at the edges $e_{1}$ and $e_{2}$, respectively, that are inclined closest to the switch edge. The dotted curve in all of these figures, which we denote by $L_{F}$, represents the intersection of the planar development $\Pi_{F}$ of $F$ with the plane that (purportedly) contains the turnover determined by $\Pi_{F}$, $\Pi_{1}$ and $\Pi_{2}$. 

\begin{remark}\label{R:ParityAmbiguity}
The symbol ``$^{*}$'' attached to a letter in this figure and in all subsequent figures is meant to indicate an ambiguity that may arise due to parity, and it is important for us to take note of it. For example, in Figure \ref{F:OneSwitch2}, if the order of the edge $AB$ is even, then the vertex $C^{*}$ is a developed copy of the vertex $C$, and the order of the edge $AC^{*}$ is also $q$, i.e., the order of edge $AC$. However, if $l$ is odd, then it would take an odd number $l$ of tetrahedra developed around the edge $AB$ to continue the development of the face $ABC$, making $C^{*}$ a developed copy of the vertex $D$ (recall that, behind the page, relative to the reader, lies the fourth vertex $D$ of the tetrahedron), and making the order of the edge  $AC^{*}$ equal to $p$, i.e., the order of the edge $AD$ (recall the notation $T[l,m,q;n,p,r]$ defined in Figure \ref{F:TLMQNPR}). We will avoid this notation whenever it is possible, although it will be necessary at times.
\end{remark}

By the previous case, we know that $\Pi_{1}$ meets none of the planes through edge $BC$. It is therefore necessary, if $\Pi_{1}$ and $\Pi_{2}$ are to intersect, that $\Pi_{2}$ cross every plane through edge $BC$. As in the previous case, then, we can conclude that one of the edges incident at $B$ must have order 2, for otherwise it is not possible for $\Pi_{2}$ to cross the plane through $BC$ inclined closest to the switch.

\begin{figure}
	\centering
	\def\svgwidth{6in}
	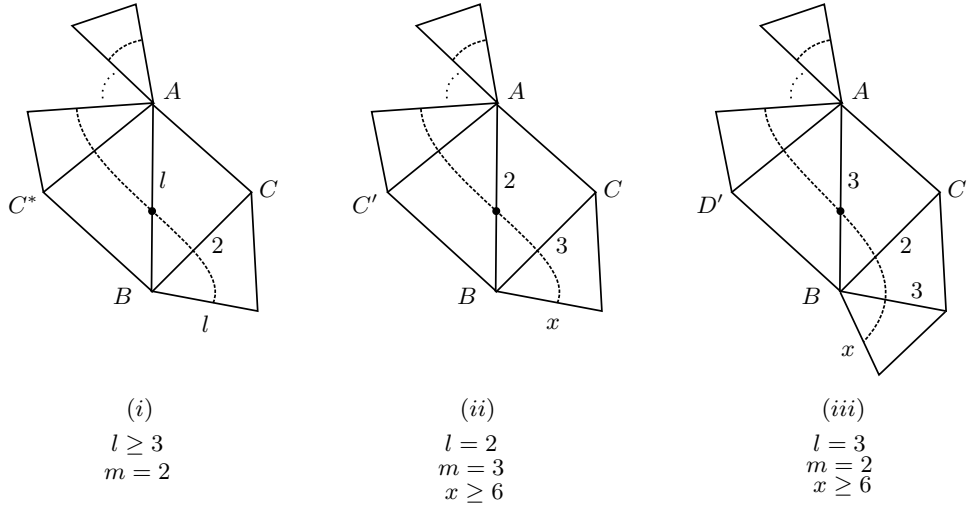
	\caption{After analysis, the remaining cases of Figure \ref{F:FaceSubdiagrams1}(c).}
	\label{F:OneSwitch2Cases}
\end{figure}

\begin{figure}
	\centering
	\def\svgwidth{5in}
	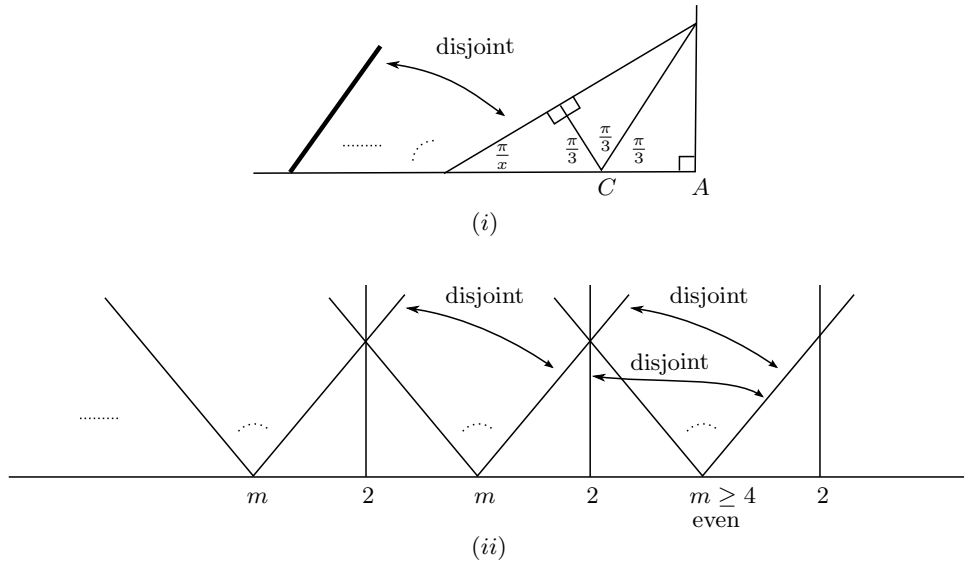
	\caption{Patterns of intersections of certain lines corresponding to sides in a triangular tiling of $\mathbb{H}^{2}$ or $\mathbb{E}^{2}$. Double arrows indicate two lines that do not intersect above the horizontal line.}
	\label{F:UseLineIntersectionsCase2}
\end{figure}

Using Figure \ref{F:TriangleLineIntersections} and the fact that $B$ must have an incident order 2 edge, we can reduce the cases that must be considered to those listed in Figure \ref{F:OneSwitch2Cases}, as follows.
Referring to Figure \ref{F:OneSwitch2}, suppose first that $l=2$ and $m=3$. Recall that the dotted curve represents the line $L_{F}$. Then the next edge incident to $B$ that  $L_{F}$ crosses after $BC$ in the direction away from the switch should have order $x \geq 6$. A schematic of the view from $B$ is pictured in Figure \ref{F:UseLineIntersectionsCase2}$(i)$. The bold line in the figure represents \emph{any} plane through a subsequent edge incident to $B$ that $L_{F}$ crosses after the edge with order $x$.  
Because the angle $\alpha$, which is formed by the bold line and the line $AC$, will always be at least $\pi/x$, we conclude that the two lines indicated in the figure by the endpoints of the double arrow will not intersect above the line $AC$. Consequently, because the line $AC$ represents the plane $\Pi_{F}$, we conclude that the planes represented by these lines will not intersect on the other side of $\Pi_{F}$ (recall that the other side of $\Pi_{F}$ refers to the side underneath the page in Figure \ref{F:OneSwitch2}). Therefore, we have reduced the case of showing that $\Pi_{1} \cap \Pi_{2} = \emptyset$ in Figure \ref{F:OneSwitch2} to the case of subfigure $(ii)$ in Figure \ref{F:OneSwitch2Cases}, provided that $l=2$ and $m=3$. The case when $l=2$ and $m$ is even with $m \geq 4$ can be eliminated in an entirely similar fashion. See Figure \ref{F:UseLineIntersectionsCase2}$(ii)$, which shows the pattern of intersections of lines that would result in the view from $B$. Here, we consider the right-most point on the horizontal (the horizontal represents $\Pi_{F}$ in the view from $B$) with the label 2 as corresponding to the edge $AB$, and the right-most point on the horizontal with the label $m$ as corresponding to the edge $BC$. It is readily seen from the figure that no lines passing through the labeled points on the horizontal to the left of the right-most point labeled $m$ ever intersect the line through the latter point that is inclined closest to the switch point (i.e., the right-most point labeled 2). Therefore, no plane through an edge incident to $B$ that is crossed by $L_{F}$ after the edge $BC$ can intersect the plane through $BC$ inclined closest to the switch, when $l=2$ and $m$ is even and at least 4. Therefore, no plane through an edge incident to $B$ that is crossed
by $L_F$ after the edge $BC$ (such as $\Pi_{2}$) can intersect the
plane $\Pi_{3}$ through $BC$ inclined closest to the switch, when $l=2$
and $m$ is even and at least 4.  Since $\Pi_{1}$ will also be disjoint
from $\Pi_{3}$ (by Subsection \ref{SSS:FaceSubdiagrams1b}), $\Pi_{1}$ will be separated from
$\Pi_{2}$ by $\Pi_{3}$, which eliminates this case.
 In fact, all of the other reductions are arrived at in this way, that is, by using the information in Figure \ref{F:TriangleLineIntersections}. The other cases that are \emph{eliminated} by the methods of this paragraph are: (1) $l=2$ and $m \geq 5$ with $m$ odd, (2) $l=3$ and $m \geq 6$ and (3) $l \geq 6$ and $m=3$. The other cases that are \emph{reduced} by the methods of this paragraph are: (4) $l \geq 3$ and $m=2$ (which reduces to the case of Figure \ref{F:OneSwitch2Cases}$(i)$) and (5) $l=3$ and $m=2$ (which reduces to the case of Figure \ref{F:OneSwitch2Cases}$(iii)$). (We note that, when $l=3$ and $m=2$, case $(i)$ of Figure \ref{F:OneSwitch2Cases} may seem to rule out case $(iii)$. However, the plane inclined closest to the switch through the edge labeled $x$ in case $(iii)$ intersects the plane inclined closest to the switch through the lower edge labeled 3 (this may be seen using the information of Figure \ref{F:TriangleLineIntersections}). We therefore must show that $\Pi_{1} \cap \Pi_{2} = \emptyset$ in \emph{both} the case that $e_{2}$ is the lower edge labeled $l=3$ in $(i)$ \emph{and} in the case that $e_{2}$ is the lower edge labeled $x$ in $(iii)$.) 

\begin{figure}
	\centering
	\def\svgwidth{5in}
	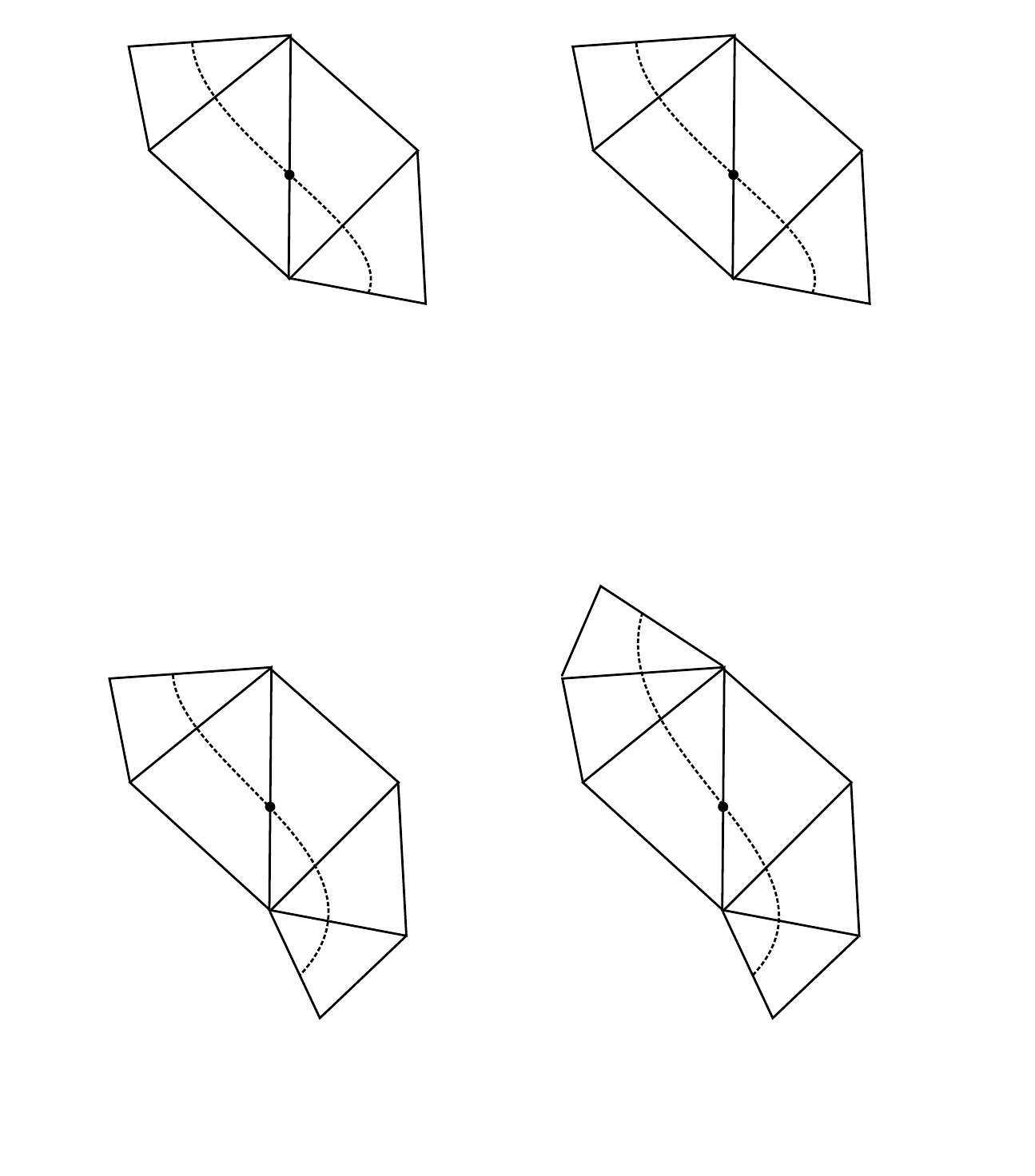
	\caption{After further analysis, applied to the cases of Figure \ref{F:OneSwitch2Cases}, these are the remaining cases of Figure \ref{F:FaceSubdiagrams1}(c) to consider.}
	\label{F:OneSwitch2CasesPart2}
\end{figure}

Now, we apply the arguments of the previous two paragraphs to the other direction along $L_{F}$ from the switch. Specifically, referring to Figure \ref{F:OneSwitch2}, we know by the previous case that $\Pi_{2}$ meets none of the planes through the edge $AC^{*}$, and so we reduce the possibilities for the number of developed faces around the vertex $A$ using the fact that $\Pi_{1}$ must intersect every plane through the edge $AC^{*}$ in order for it to be possible for $\Pi_{1}$ and $\Pi_{2}$ to have nonempty intersection. The result of this further analysis leaves us to consider only the cases of Figure \ref{F:OneSwitch2CasesPart2}. We note the change from ``$l \geq 3$'' to ``$l \geq 3$ odd'' that occurs when reducing Figure \ref{F:OneSwitch2Cases}$(i)$ to Figure \ref{F:OneSwitch2CasesPart2}$(i)$. This change is due to the fact that, when $l$ is even, the edge label ``2'' for $AD'$ in \ref{F:OneSwitch2CasesPart2}$(i)$ must equal the edge label for $AC$. However, this would contradict our assumption that none of the vertices of $T$ is finite, because $C$ would have two incident edges, $AC$ and $BC$, labeled 2. 

\begin{figure}
	\centering
	\def\svgwidth{5in}
	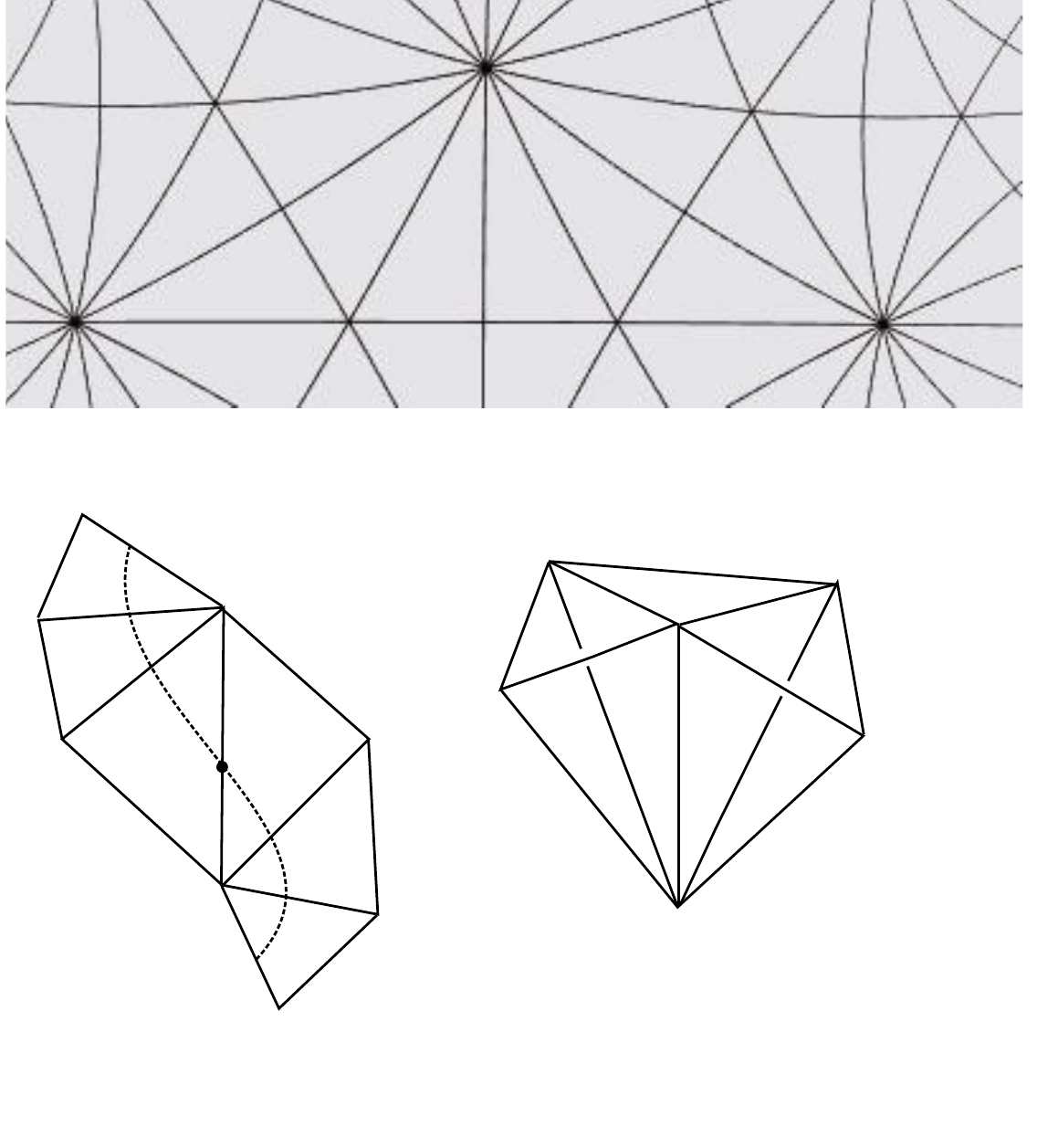
	\caption{The case of Figure \ref{F:OneSwitch2CasesPart2}$(iv)$. The upper half of the figure represents the view from the vertex $A$ when $y=7$. The lower half consists of a perspective image of the three copies of the tetrahedron $ABCD$ on the right, and several triangles in the development of the face $ABC$ on the left.}
	\label{F:OneSwitch2CasesPart2_CaseIV}
\end{figure}

So we are left to analyze the cases of Figure \ref{F:OneSwitch2CasesPart2}. We begin with case $(iv)$. See Figure \ref{F:OneSwitch2CasesPart2_CaseIV}. The multiple parts of this figure are explained in the caption.
Referring to the left side of the lower half of the figure, $\Pi_{1}$ is the plane through edge $AC''$ inclined closest to the switch edge $AB$ and $\Pi_{2}$ is the plane through edge $BD''$ inclined closest to the switch edge $AB$. We wish to show that $\Pi_{1} \cap \Pi_{2} = \emptyset$. We do so using the upper half of the figure, which shows the view from $A$ under the assumption that $y=7$ (the same argument we give here applies to any other value for $y \geq 6$). In the upper half of the figure, the plane $\Pi_{1}$ is represented by the line $C''D'''$, and the plane $ACD$---which is depicted in the right side of the lower half of the figure, and which is the plane through $AC$ inclined closest to the switch edge $AB$ in the left side of the lower half of the figure---is represented by the line $CD$. Recalling that $\Pi_{F}$ is the plane containing the face $ABC$ (and, therefore, the plane in which the left side of the lower half of the figure is drawn, as well as the horizontal line in the upper half of the figure), we observe that there are two planes, other than $\Pi_{F}$, that pass through $AB$. These planes are represented in the upper half of the figure by the lines $BC'$ and $BD$. Using the upper half of the figure, we observe that any point of $\Pi_{1}$ that is on the same side of $ACD$ as the vertex $B$ is also on the same side of the plane $ABC'$ (which is represented by the line $BC'$) as the point $D'$. We now use the previous case (Section \ref{SSS:FaceSubdiagrams1b}) to observe that $\Pi_{2} \cap ACD = \emptyset$: namely, $\Pi_{2}$ and $ACD$ are the planes through $BD''$ and $AC$, respectively, inclined closest to the \emph{new} switch edge $BC$ for the three triangles $ABC$, $A'BC$ and $A'BD''$ from the lower left half of Figure \ref{F:OneSwitch2CasesPart2_CaseIV}, to which Section \ref{SSS:FaceSubdiagrams1b} applies (to see this more clearly, turn these three triangles together so that the edge $BC$ is vertical, and compare with Figure \ref{F:OneSwitch1}). In exactly the same way (i.e., using Section \ref{SSS:FaceSubdiagrams1b}), we see that $\Pi_{2} \cap AC'D' = \emptyset$, this time using $AB$ as the switch edge. But now, since $\Pi_{2}$ is on the same side of $ACD$ as the vertex $B$ \emph{and} on the same side of $AC'D'$ as the vertex $B$, we can use the upper half of Figure \ref{F:OneSwitch2CasesPart2_CaseIV} to see that there is no part of $\Pi_{1}$ which is both on the $B$ side of $AC'D'$
\emph{and} on the $B$ side of $ACD$. Therefore, $\Pi_{1} \cap \Pi_{2} = \emptyset$.

\begin{figure}
	\centering
	\def\svgwidth{6.5in}
	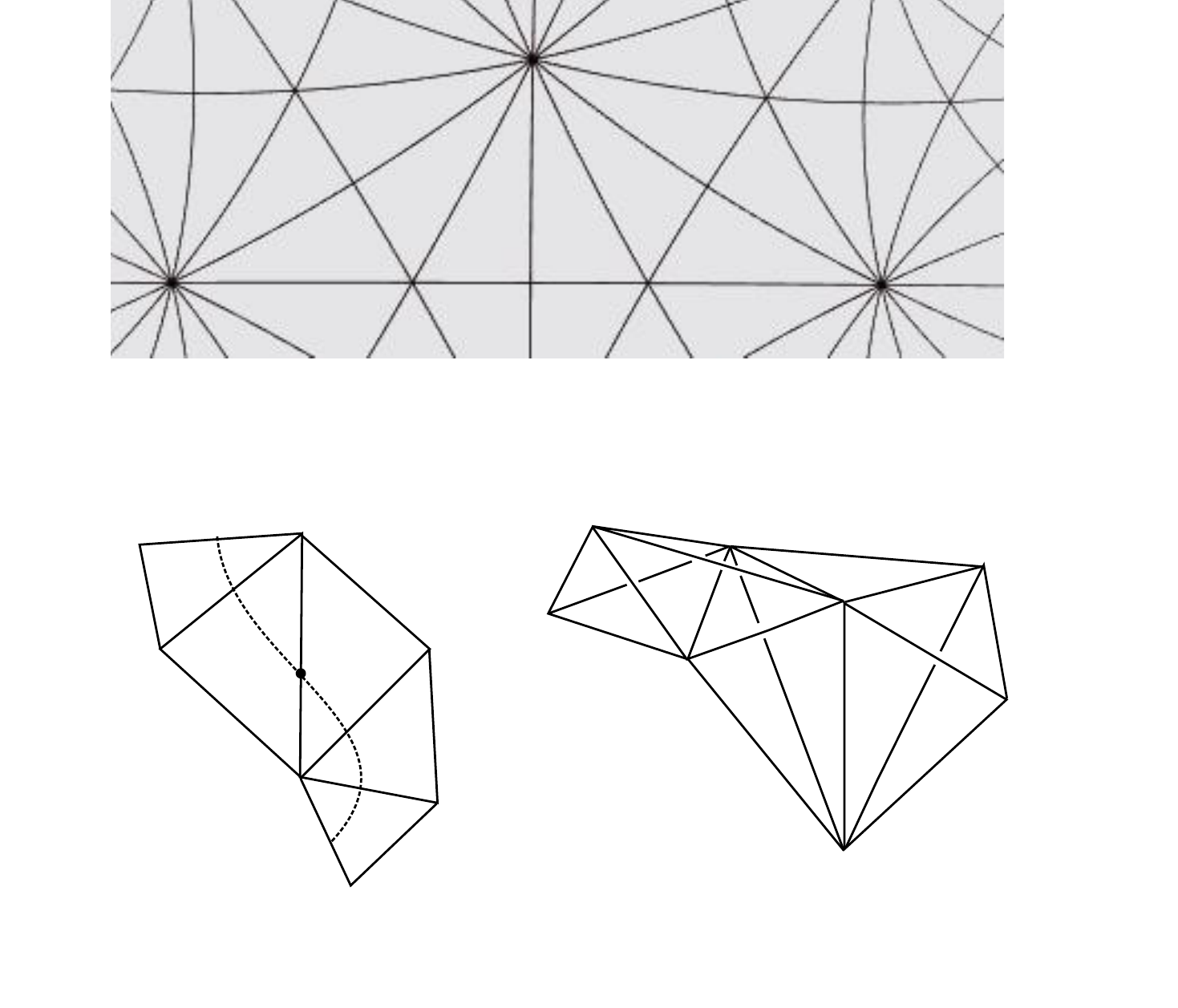
	\caption{The case of Figure \ref{F:OneSwitch2CasesPart2}$(iii)$. The upper half of the figure represents the view from vertex $A$ when $y=7$. The right side of the lower half of the figure depicts the development of several copies of the tetrahedron, and the left side of the lower half depicts the development of the face $ABC$.}
	\label{F:OneSwitch2CasesPart2_CaseIII}
\end{figure}

The argument of the previous paragraph can be used in case $(iii)$ of Figure \ref{F:OneSwitch2CasesPart2}. See Figure \ref{F:OneSwitch2CasesPart2_CaseIII}. In the lower left half of this figure, $\Pi_{1}$ is the plane through the edge $AB'$ inclined closest to the switch edge $AB$. In the lower right half, $\Pi_{1}$ is the plane $AC'B'A''$. In the upper half of the figure, which represents the view from $A$ when $y=7$ (the case when $y \geq 6$ is similar), $\Pi_{1}$ is represented as the line $B'C'$. 
Proceeding as in the previous paragraph, we have $\Pi_{2} \cap ACD = \Pi_{2} \cap AC'D' = \emptyset$ (by Section \ref{SSS:FaceSubdiagrams1b}). Furthermore, $\Pi_{2}$ is on the $B$ side of both $ACD$ and $AC'D'$. But now, referring to the upper half of Figure \ref{F:OneSwitch2CasesPart2_CaseIII}, we see that there is no part of $\Pi_{1}$ that is on the $B$ side of both $ACD$ and $AC'D'$. So $\Pi_{1} \cap \Pi_{2} = \emptyset$.

\begin{figure}
	\centering
	\def\svgwidth{5in}
	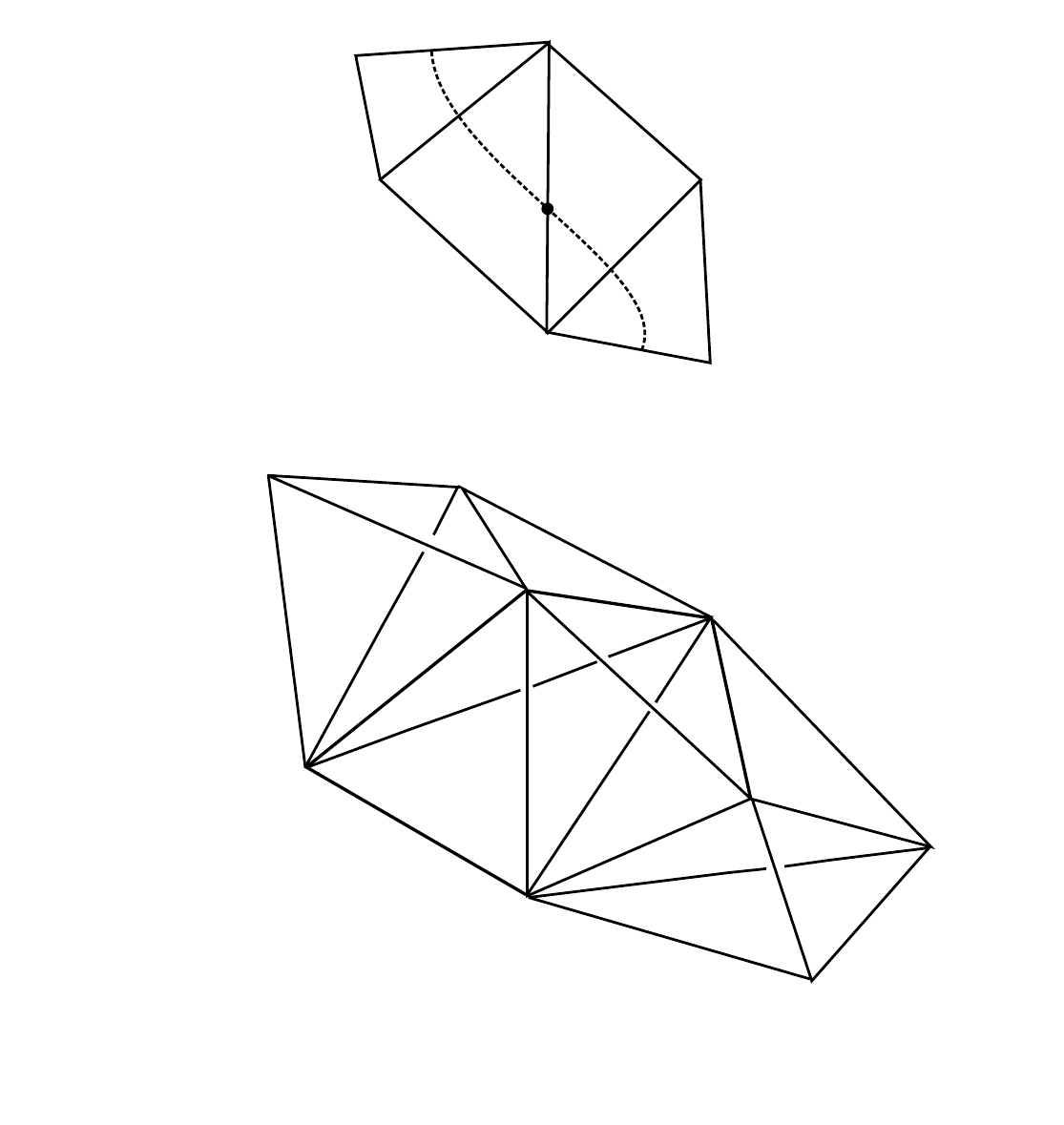
	\caption{The case of Figure \ref{F:OneSwitch2CasesPart2}$(ii)$.}
	\label{F:OneSwitch2CasesPart2_CaseII_REDUX}
\end{figure}

\begin{figure}
	\centering
	\def\svgwidth{4.5in}
	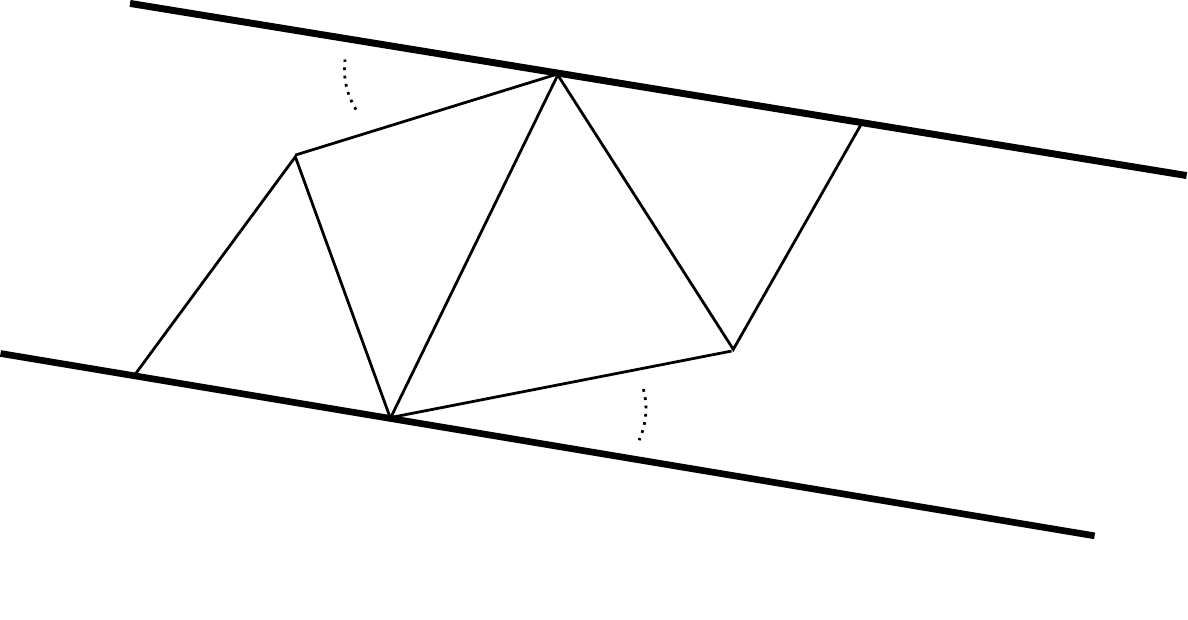
	\caption{The schematic view from the vertex $D$ for Figure \ref{F:OneSwitch2CasesPart2_CaseII_REDUX}. The letters $x$, $y$ and $n$ represent the integral submultiples of $\pi$ of the dihedral angles of the tetrahedra incident at $D$.}
	\label{F:OneSwitch2CasesPart2_CaseII_REDUX2}
\end{figure}

We now address case $(ii)$ of Figure \ref{F:OneSwitch2CasesPart2}. See Figure \ref{F:OneSwitch2CasesPart2_CaseII_REDUX}. 
Referring to the upper part of this figure, $\Pi_{1}$ and $\Pi_{2}$ are the planes through the edges $AD'$ and $BD''$, respectively, that are inclined closest to the switch edge $AB$. In the lower part of the figure, which depicts the development of multiple copies of the tetrahedron, $\Pi_{1}$ is the plane $ADB'D'$ and $\Pi_{2}$ is the plane $BDA'D''$. Because these two planes are both incident to the non-finite vertex $D$, they intersect if and only if their intersections with the link of $D$ intersect. See Figure \ref{F:OneSwitch2CasesPart2_CaseII_REDUX2}, which schematically depicts the view of this link from the vertex $D$. 
In the figure, $\Pi_{1}$ is represented by the bold line $AB'$ and $\Pi_{2}$ by the bold line $A'B$. We have labeled the interior angles of the triangles in this view by their submultiples of $\pi$. Because vertex $C$ has two edges of order 3 incident to it, we must have that $n\geq3$. But since $x$ and $y$ must both be at least $6$, we can use Figure \ref{F:TriangleLineIntersections}$(iii)$ (with base the segment $AB$) to conclude that the bold lines cannot intersect on either side of the line $AB$. So $\Pi_{1} \cap \Pi_{2} = \emptyset$ in the case of Figure \ref{F:OneSwitch2CasesPart2}$(ii)$.    

\begin{figure}
	\centering
	\def\svgwidth{6.5in}
	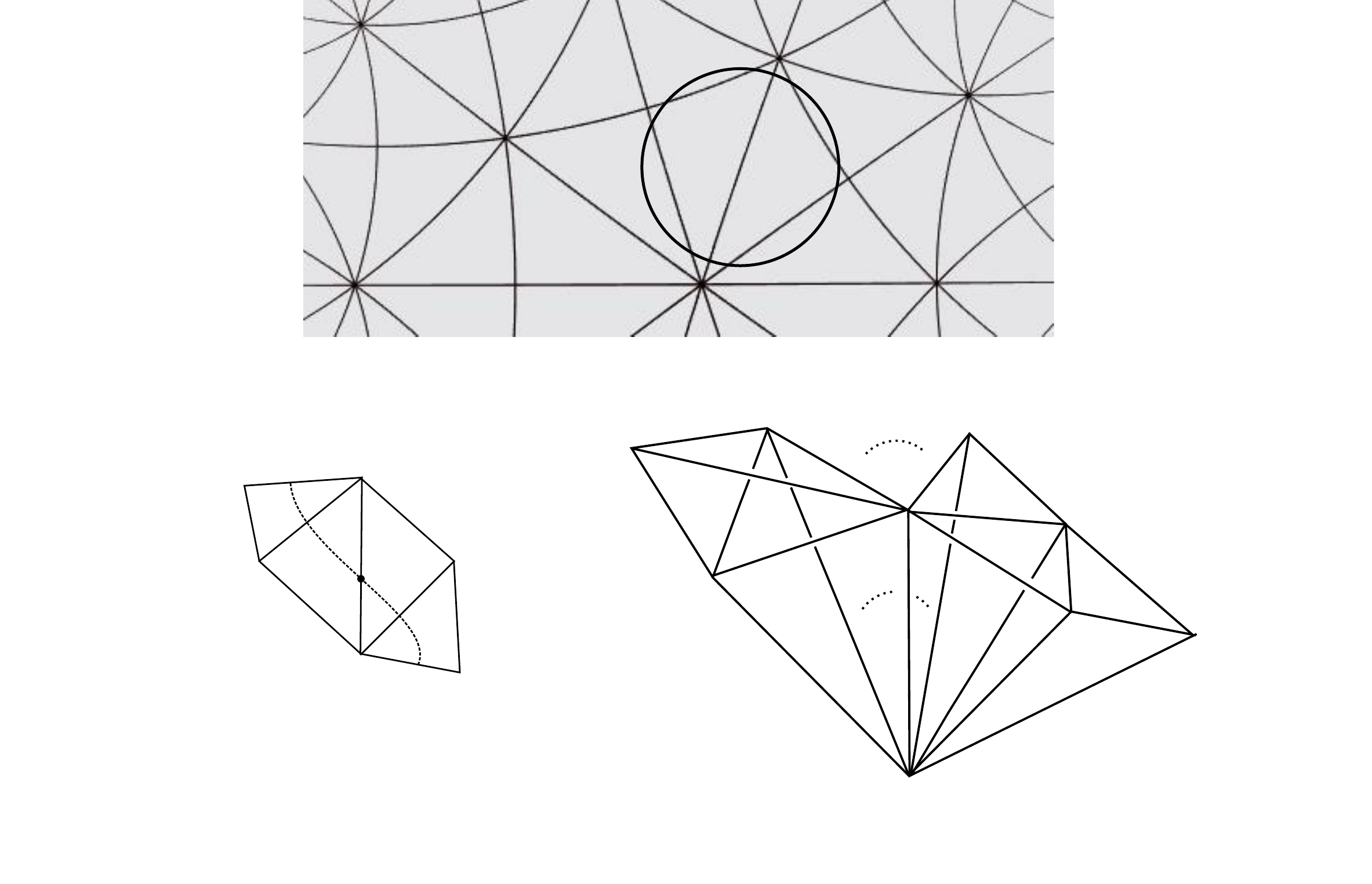
	\caption{The case of Figure \ref{F:OneSwitch2CasesPart2}$(i)$. The upper half of the figure depicts the view from the vertex $A$, in the case when $A$ has type $(2,4,5)$.}
	\label{F:OneSwitch2CasesPart2_CaseI}
\end{figure}

This leaves case $(i)$ of Figure \ref{F:OneSwitch2CasesPart2}. We begin by assuming that $l \geq 5$ (recall that $l$ must be odd). See Figure \ref{F:OneSwitch2CasesPart2_CaseI}. 
The upper half of this figure depicts the view from the vertex $A$ with the projection centered at the vertex $B$. For the purposes of illustration, we take the type of $A$ to be $(2,4,5)$, although the argument only depends on the presence of the order 2 edge incident to $A$ and the fact that the order of the edge $AB$ is at least 5. The plane $\Pi_{1}$ is represented by the circular arc $B'C''$. We note that the plane $BC'D$ in the lower right part of the figure is represented in the upper half of the figure by a circle, centered on the line segment $BC'$ because the planes $BC'D$ and $ABC'$ are perpendicular, and whose interior disk does not contain any of the points $B$, $D$, $C'$ or $D'$. As we have observed previously (see the argument depicted in Figure \ref{F:RuleOut2xy_1} from Section \ref{SSS:SingleEdgeCrossedOrder4OrMore}), the circle representing $BC'D$ can intersect at most two sides of the $l$--gon centered at $B$, and in this case those sides will always be $D'C'$ and $DC'$. It is clear that this circle is disjoint from the arc $B'C''$ representing $\Pi_{1}$, and hence that $\Pi_{1} \cap BC'D = \emptyset$. Now referring to the lower right part of the figure, we observe that $\Pi_{2}=A'BD$ (as planes) and that the part of $\Pi_{2}$ that is on the same side of $BCD$ as $A$ is also on the \emph{opposite} side of $BC'D$ as $A$. Since $\Pi_{1}$ is disjoint both from $BC'D$ and $BCD$ (the latter by the previous case of Section \ref{SSS:FaceSubdiagrams1b}), and because $\Pi_{1}$ lies on the same side of these planes as $A$, we can conclude that $\Pi_{1} \cap \Pi_{2} = \emptyset$ in this case when $l \geq 5$.

\begin{figure}
	\centering
	\def\svgwidth{6.5in}
	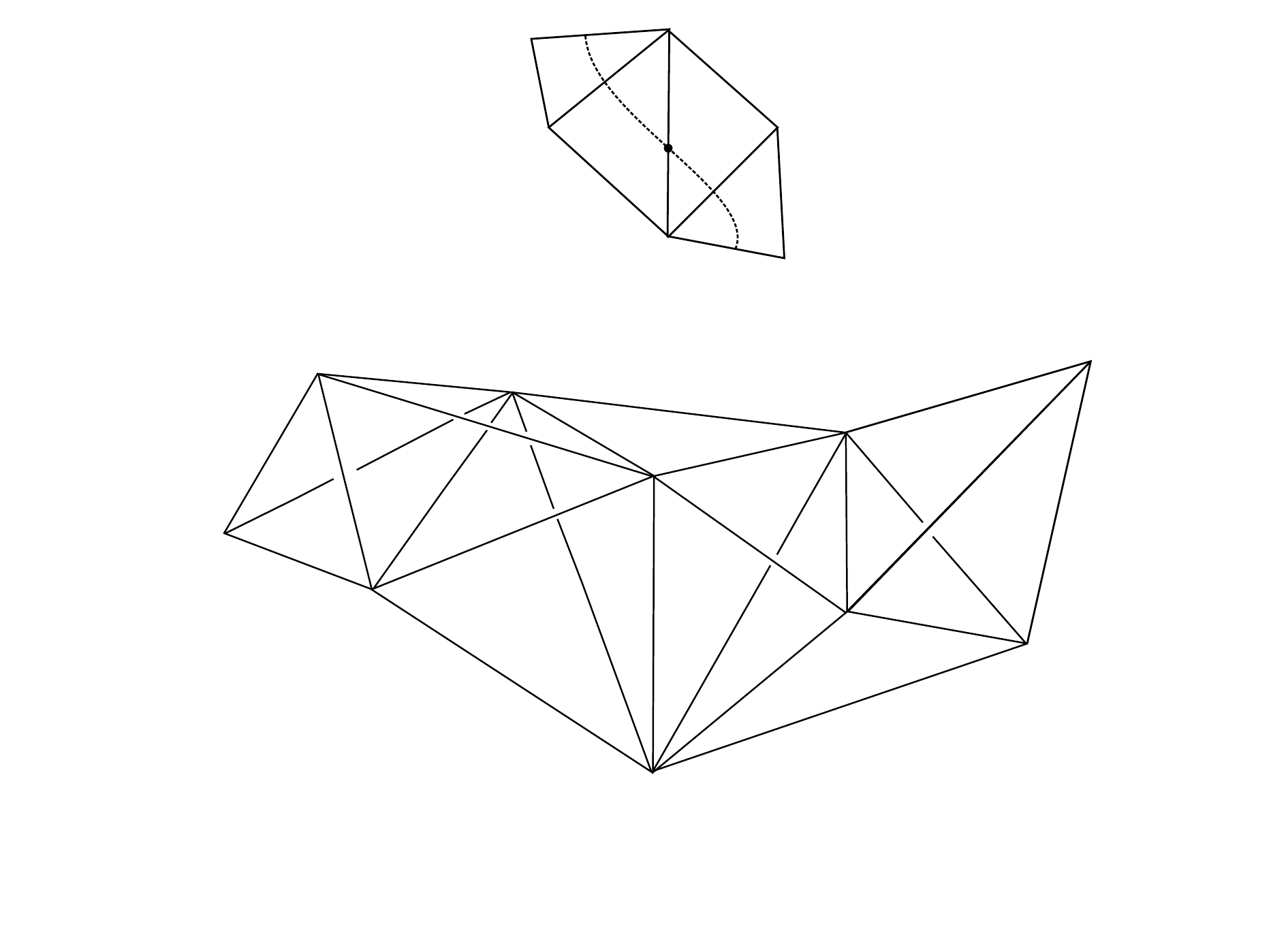
	\caption{The case of Figure \ref{F:OneSwitch2CasesPart2}$(i)$ when $l=3$.}
	\label{F:OneSwitch2CasesPart2_CaseI2}
\end{figure}

So we now assume that $l=3$ in this case. We are not able to use the argument of the previous paragraph because some of the intersections ruled out in the previous paragraph can occur in this case. We refer to Figure \ref{F:OneSwitch2CasesPart2_CaseI2}.
The possible values for $q$, $n$ and $r$ in the figure are based on the fact that the tetrahedron has no finite vertices. In this figure, $\Pi_{1} = AC'A''B'$ and $\Pi_{2}=A'B''DB$ (as planes). We determine that these planes are disjoint by applying the techniques of Sections \ref{SSS:SingleEdgeCrossedOrder2} and \ref{SSS:SingleEdgeCrossedOrder3}. In particular, if $n \geq 4$, then we use the geometry of the link of vertex $C'$ to conclude that $\Pi_{1}$ is disjoint from the plane $BDC'D'$ (it lies to the same side of $BDC'D'$ as the vertex $A$) and the geometry of the link of vertex $D$ to conclude that $\Pi_{2}$ is disjoint from the plane $ACDC'$ (it lies to the same side of $ACDC'$ as vertex $B$). Now by considering the plane $ABD$ and the geometry of the vertex $B$, we have that the part of $\Pi_{2}$ that is on the $C'$ side of $ABD$ is always on the opposite side of $BDC'D'$ to $\Pi_{1}$. Similarly, we have that the part of $\Pi_{1}$ on the $D$ side of $ABC'$ is always on the opposite side of $ACDC'$ to $\Pi_{2}$. We conclude that $\Pi_{1} \cap \Pi_{2} = \emptyset$. When $n=3$, the argument is similar, except that $ACDC' = ACB''DC'$ and $BDC'D' = BDC'A''D'$ (as planes), and  $\Pi_{1}$ and $\Pi_{2}$ will  form interior angles on the $B$ side of $ACB''DC'$ of $3\pi/q \leq \pi/2$ and $\pi/r \leq \pi/6$, respectively (so that $\Pi_{1}$ and $\Pi_{2}$ cannot intersect on the side of this plane opposite to $B$), and interior angles on the $A$ side of $BDC'A''D'$ of $\pi/q \leq \pi/6$ and $3\pi/r \leq \pi/2$, respectively (so that $\Pi_{1}$ and $\Pi_{2}$ cannot intersection on the side of this plane opposite to $A$). Again, we conclude that $\Pi_{1} \cap \Pi_{2} = \emptyset$. This completes case $(i)$ of Figure \ref{F:OneSwitch2CasesPart2}, and concludes this subsection.      

\subsubsection{\ref{F:FaceSubdiagrams2}(a):}\label{SSS:FaceSubdiagrams2a} 
\begin{figure}
	\def\svgwidth{3in}
	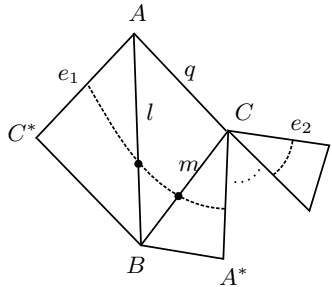
	\caption{One case of Figure \ref{F:FaceSubdiagrams2}(a).}
	\label{F:TwoSwitch1}
\end{figure}

See Figure \ref{F:TwoSwitch1}, and recall the significance of the symbol ``$^{*}$'' from Remark \ref{R:ParityAmbiguity}. We must first address the case when $e_{2} = A^{*}C$. 
There are two possibilities that we must consider in determining whether or not $\Pi_{1}$ and $\Pi_{2}$ can intersect: either (1) $\Pi_{1}$ meets the plane through $BC$ that is closest in inclination to the switch edge $AB$ (it cannot meet more planes through $BC$, by our previous observations) and $\Pi_{2}$ meets at least the second closest plane through $BC$ to the switch edge $AB$, or (2) $\Pi_{1}$ meets no planes passing through $BC$ and $\Pi_{2}$ meets all of the planes passing through $BC$. We handle these two cases below:
\begin{enumerate}
	\item In order for $\Pi_{1}$ to meet a plane passing through $BC$, our tetrahedron must take one of the forms of items (1)--(3) in the summary at the conclusion of the paper. This follows from the extensive analysis of Section \ref{SS:SingleEdgeCrossed} (in fact, the pairwise intersections of $\Pi_{1}$, $\Pi_{F}$ and the plane through $BC$ inclined closest to the switch edge $AB$ determine an immersed turnover in this case).  We consider the case when $l=3$, corresponding to item (3) in the summary. If $l=3$, then $q=2$, $m \geq 6$ and $n$ (the order of the third edge associated to vertex $C$) is at least 3. It is then an easy analysis, using Figure \ref{F:TriangleLineIntersections} applied to the vertex $C$, to see that there is no choice of $n$ and $m$ for which $\Pi_{2}$ can intersect either of the two closest planes through $BC$ toward the edge $AB$. So $\Pi_{1} \cap \Pi_{2} = \emptyset$. Exactly the same analysis holds if our tetrahedron takes the form of item (1) of the summary at the conclusion of the paper (in this case we have $l=2$, $m \geq 6$, $n=2$ and $q \geq 3$, and so the order of edge $A^{*}C$ is either $2$ or $q$, and there is no choice for $m$ and $q$ such that $\Pi_{2}$ meets either of the two planes through $BC$ inclined closest to the switch). If our tetrahedron has the form of item (2) from the summary, then $l=2$, $m \geq 3$ and $q \geq 6$. If $m$ is odd and at least 5, then the order of edge $A^{*}C$ is $2$ and we can use Figure \ref{F:TriangleLineIntersections} applied to vertex $C$ to conclude that $\Pi_{2}$ does not meet the two planes through $BC$ inclined closest to the switch. If $m$ is even and at least 6, then the order of $A^{*}C$ is $q \geq 6$, and the conclusion of the previous sentence also holds. If $m = 4$, then we refer to Figure \ref{F:TwoSwitch1Case1}.
\begin{figure}
	\def\svgwidth{4.5in}
	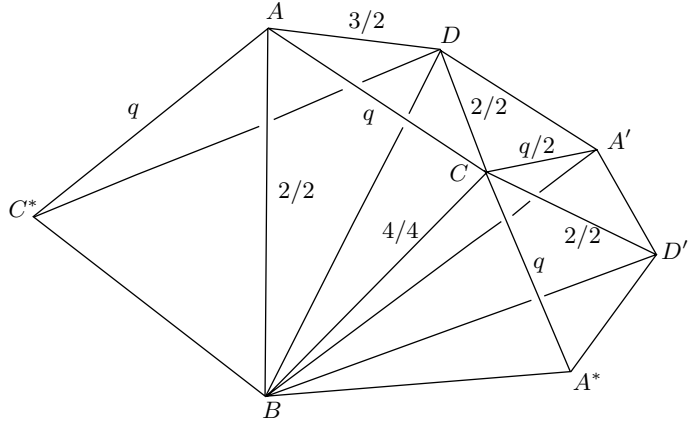
	\caption{The case of Figure \ref{F:TwoSwitch1}, when $l=2$, $m=4$ and $e_{2} = A^{*}C$.}
	\label{F:TwoSwitch1Case1}
\end{figure} 
Only the relevant edges are labeled in this figure, in which $\Pi_{1} = AC^{*}D$ and $\Pi_{2} = A^{*}D'A'C$. Because $q \geq 6$, we have that $\Pi_{1}$ and $\Pi_{2}$ form interior angles on the side of $ACA'D$ opposite to vertex $B$ of $\pi/3$ and $(q-2)\pi/q \geq 2\pi/3$, respectively (these are interior angles with respect to the edge $CD$). Therefore, $\Pi_{1}$ and $\Pi_{2}$ do not intersect on the side of this plane opposite to $B$. But, as we have observed, $\Pi_{1} \cap A'BC = \emptyset$. Since the part of  $\Pi_{2}$ that is on the $B$ side of $ACA'D$ is always on the opposite side of $A'BC$ to $\Pi_{1}$, we have $\Pi_{1} \cap \Pi_{2} = \emptyset$. The case when $m=3$ is exactly the same. These are all the possibilities for when the tetrahedron has one of the types (1)--(3) in the summary. So $\Pi_{1} \cap \Pi_{2} = \emptyset$ for this case.
	\item If the order of edge $BC$ is greater than 4, then it is not possible to choose integers for the type of vertex $C$ so that $\Pi_{2}$ crosses all the planes through $BC$. This follows by using the information of Figure \ref{F:TriangleLineIntersections} applied to the vertex $C$, as in the arguments that accompany Figure \ref{F:UseLineIntersectionsCase1} in Section \ref{SSS:FaceSubdiagrams1b}. The same statement is true (with the same argument) if the order of $BC$ is 3 and the vertex $C$ has no incident order 2 edge. So the order of edge $BC$ is either 3 and $C$ has the type $(2,3,x\geq6)$ or the order of edge $BC$ is 2. Suppose that the edge  $BC$ has order 3. Then we can use the same argument as the one given at the end of the previous paragraph. Namely, it is readily shown that $\Pi_{1}$ and $\Pi_{2}$ meet the plane containing the face $ACD$ at interior angles that sum to at least $\pi$ on the opposite side of $ACD$ of the vertex $B$, and since they do not meet on the $B$ side of this plane, they must be disjoint. The same argument also works when the order of $BC$ is 2. So $\Pi_{1} \cap \Pi_{2} = \emptyset$ in this case.
\end{enumerate}

So we assume $e_{2} \neq A^{*}C$. We observe that removing the sides $AC^{*}$ and $BC^{*}$ from the Figure \ref{F:TwoSwitch1} leaves a picture that is equivalent to the previous case of Subsection \ref{SSS:FaceSubdiagrams1b}. We therefore know that $\Pi_{2}$ misses every plane through the switch edge $AB$. It follows, using Figure \ref{F:TriangleLineIntersections} applied to the vertex $A$, that $l$ must be either 2 or 3, in order for $\Pi_{1}$ to cross every plane through this switch edge. Moreover, we must have, as in previous cases, that the type of vertices $A$ and $C$ must include an order 2 point.  Suppose $l=2$. This implies that neither $m$ nor $q$ is 2. If, in addition, neither $m$ nor $q$ is 3, then it is straightforward using the information in Figure \ref{F:TriangleLineIntersections} (applied to vertex $C$) to show that $\Pi_{2}$ cannot meet the plane through $A^{*}C$ inclined closest to $\Pi_{1}$, and so prove that $\Pi_{1} \cap \Pi_{2} = \emptyset$ in this case. So either $m=3$ and $q \geq 6$ or $q=3$ and $m \geq 6$, and in both cases $n=2$. In either case, it is a straightforward application of the techniques already employed---specifically, the techniques involving developing tetrahedra from Sections \ref{SSS:SingleEdgeCrossedOrder2} and \ref{SSS:SingleEdgeCrossedOrder3}---to show that $\Pi_{1}$ and $\Pi_{2}$ do not intersect.

\begin{figure}
	\def\svgwidth{4.5in}
	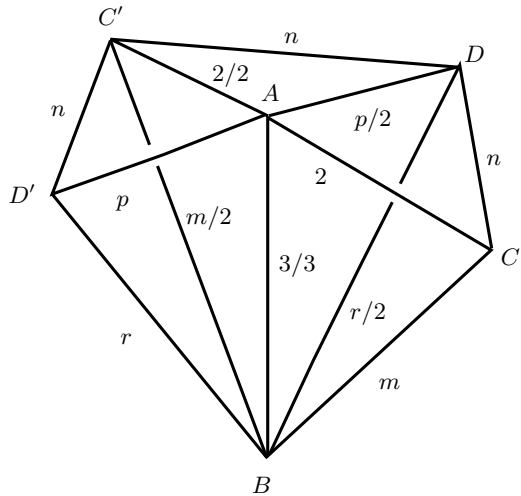
	\caption{The case of Figure \ref{F:TwoSwitch1}, when $l=3$, $q=2$ and $p\geq6$.}
	\label{F:TwoSwitch1Case2}
\end{figure}

Now suppose $l=3$. This implies that the order of edge $e_{1}$ is $p$. Because $\Pi_{1}$ must cross every plane through the switch edge $AB$, it is easily shown using Figure \ref{F:TriangleLineIntersections} (applied to vertex $A$) that the order $p$ of edge $e_{1}$ is at least 6 and $q = 2$. Figure \ref{F:TwoSwitch1Case2} shows three copies of $T$, with $C^{*}$ relabeled as $D'$. Because $q=2$, we must have $n\geq3$ and $m\geq3$. 
Since $n \neq 2$, analysis using the vertex $D$ shows that $\Pi_{1}$, which is the plane $ADC'D'$, intersects the plane $BCD$ if and only if $r=2$. We analyze the two cases $r\neq2$ and $r=2$:
\begin{enumerate}
	\item $r \neq 2$:  In this case, $\Pi_{1}$ does not intersect $BCD$, and so it is necessary for $\Pi_{2}$ to intersect $BCD$ if $\Pi_{1}$ and $\Pi_{2}$ are to intersect. If $m \geq 4$ and even, then the edges emanating from the vertex $C$ in Figure \ref{F:TwoSwitch1}---$CA$, $CB$, $CA^{*}$,..., $e_{2}$---have labels that alternate $2,m,2,...$ However, by using Figure \ref{F:UseLineIntersectionsCase2}$(ii)$ applied to the vertex $C$, it is easily seen that no plane through any of the edges $CA^{*}$,..., $e_{2}$ that is inclined closest to the switch edge $AB$ will intersect the plane $BCD$. Since $\Pi_{1}$ does not intersect $BCD$, the latter plane separates $\Pi_{1}$ from $\Pi_{2}$. So we are left to consider when $m\geq3$ and odd. When $m\geq5$ and odd, an application of the information from Figure \ref{F:TriangleLineIntersections} to the vertex $C$ shows that no plane that is inclined closest to the switch edge $AB$ through any of the edges from $CA^{*}$ to $e_{2}$ can intersect with the plane $BCD$. So again, $\Pi_{1} \cap \Pi_{2} = \emptyset$. Finally, when $m=3$, it is necessary for $n$ (the label of the third edge of $T$ that meets the vertex $C$, and the label of the edge $CA^{*}$) to be at least $6$. So the type of the vertex $C$ is $(2,3,n\geq6)$, and no plane through any edge after $CA^{*}$ and up to and including $e_{2}$ that is inclined closest to the switch edge $AB$ will intersect the closest such inclined plane through the edge  $CA^{*}$ (as in Figure \ref{F:UseLineIntersectionsCase2}$(i)$). Since, by the observation of the first paragraph of this section, the closest inclined plane to the switch edge $AB$ through $CA^{*}$ is disjoint from $\Pi_{1}$, we again have $\Pi_{1} \cap \Pi_{2} = \emptyset$. This completes the analysis of the case when $r \neq 2$.
	\item $r=2$: In this case, $\Pi_{1}$ does intersect the plane $BCD$. Because $r=2$ and $l=3$, it is necessary that $m \geq 6$. We have previously observed that $\Pi_{1}$ cannot intersect with the second-closest inclined plane to the switch edge $AB$ through $BC$ (because the planes $\Pi_{1}$, $ABC$ and $BCD$ form pairwise angles of intersection $\pi/p$, $\pi/m$ and $\pi/n$, with $p\geq6$, $m\geq6$ and $n\geq3$). However, vertex $C$ has type $(2,m\geq6,n\geq3)$, and it is easily seen using the information of Figure \ref{F:TriangleLineIntersections} applied to $C$ that no plane that is inclined closest to the switch edge $AB$ through any of the edges from $CA^{*}$ to $e_{2}$ can intersect the second-closest inclined plane to $AB$ through $CB$, \emph{provided} that $m\geq7$. So this second-closest inclined plane through $CB$ separates $\Pi_{1}$ from $\Pi_{2}$, when $m \geq 7$. This leaves the case when $m=6$. But this case is handled by an argument similar to the accompanying argument for Figure \ref{F:RuleOut23x_2} in Section \ref{SSS:FaceSubdiagrams1b}. This completes the case when $r=2$, and concludes this subsection.
\end{enumerate}

\subsubsection{\ref{F:FaceSubdiagrams2}(b):}\label{SSS:FaceSubdiagrams2b} 
\begin{figure}
	\centering
	\def\svgwidth{3in}
	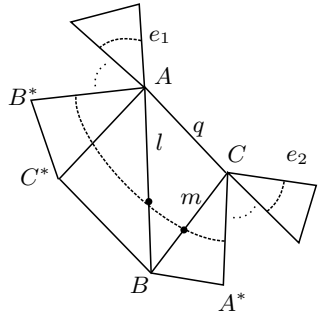
	\caption{One case of Figure \ref{F:FaceSubdiagrams2}(b).}
	\label{F:TwoSwitch2}
\end{figure}

See Figure \ref{F:TwoSwitch2}. By the result of Section \ref{SSS:FaceSubdiagrams1b}, it is not possible for $e_{2}$ to equal $CA^{*}$. Because of this, it is not possible, also by the Section \ref{SSS:FaceSubdiagrams1b}, for $\Pi_{2}$ to meet any of the planes through the edge $AB$. Nor is it possible, by Section \ref{SSS:FaceSubdiagrams1b}, for $\Pi_{1}$ to meet any of the planes through the edge $BC$. Consequently, the intersection of $\Pi_{1}$ and $\Pi_{2}$ can only occur if $\Pi_{1}$ crosses every plane through $AB$ and $\Pi_{2}$ crosses every plane through $BC$. The subsequent possibilities and arguments to rule them out are all straightforward to carry out, using the techniques we have employed to this point. 
This completes the proof. \hfill \fbox{\ref{T:TurnoversInPolyhedra}}\\

\noindent \emph{Summary.} We provide a summary of the classification of immersed turnovers in the orbifold $\cal{O}_{T}$ associated to the generalized tetrahedron $T[l,m,q;n,p,r]$. These are listed in the order in which they appear in the proof, but isometric cases are indicated (the 24 isometric cases are determined by applying an element of the symmetric group $S_{4}$: any element of the symmetric group $S_{3}$ may be applied to both the first and second triples of $T[l,m,q;n,p,r]$, and any pair from one triple may be swapped with the corresponding pair of the other triple). We also include a conjectural list of all the immersed turnovers in hyperbolic tetrahedral orbifolds. All of these can be confirmed using the techniques of this paper, and while the author believes this list to be exhaustive, the necessary computations to determine the complete classification are somewhat extensive.
\begin{enumerate}
	\item $T[2,m,q;2,p,3]$. $\cal{O}_{T}$ contains an immersed $(q,m,p)$ turnover, where $q \geq 3, m \geq 6$ and $p \geq 6$.
	\item $T[2,m,q;2,3,r]$ (isometric to item (1)). $\cal{O}_{T}$ contains an immersed $(q,m,r)$ turnover, where $q \geq 6, m \geq 3$ and $r \geq 6$.
	\item $T[3,m,2;n,p,2]$ (isometric to item (1)). $\cal{O}_{T}$ contains an immersed $(m,n,p)$ turnover, where $m \geq 6, n \geq 3$ and $p \geq 6$.
\end{enumerate}
\emph{Conjectural list of all immersed turnovers in hyperbolic tetrahedral orbifolds:}
\begin{enumerate}
	\setcounter{enumi}{3}
	\item $T[2,m,q;2,p,3]$. $\cal{O}_{T}$ contains an immersed $(q,m,p)$ turnover for any of the following values:
	\begin{enumerate}
		\item $q=2$, $m=4$ and $p \geq 5$. In this case, $\cal{O}_{T}$ also contains
			\begin{enumerate}
			\item a $(2,p,p)$ turnover,
			\item a $(4,4,5)$ turnover if $p=5$, and
			\item a $(p/2,p,p)$ turnover if $p$ is even. 
			\end{enumerate}
		\item $q=2$, $p=4$ and $m \geq 5$ (isometric to item (4), with the same set of additional non-maximal turnovers).
		\item $q=2$, $m \geq 5$ and $p \geq 5$. In this case, $\cal{O}_{T}$ also contains
			\begin{enumerate}
			\item a $(m,m,p/2)$ turnover if $p$ is even, or
			\item a $(m/2,p,p)$ turnover if $m$ is even.
			\end{enumerate}
		\item $q,$ $m$ and $p$ are all greater than 2, and at least one is greater than 3. In this case, if two of the values are 3, then $\cal{O}_{T}$ also contains a $(x,x,x)$ turnover, where $x$ is the integer that is greater than 3.
	\end{enumerate} 
	\item $T[3,2,2;2,p,3]$. $\cal{O}_{T}$ contains an immersed $(2,p,p)$ turnover, where $p \geq 5$.
	\item $T[3,m,2;2,p,3]$. $\cal{O}_{T}$ contains an immersed $(m,p,p)$ turnover, where $m \geq 3$ and $p \geq 4$.
	\item $T[3,m,3;2,3,2]$. $\cal{O}_{T}$ contains an immersed $(3,m,m)$ turnover, where $m \geq 4$.
	\item $T[4,3,q;2,2,2]$. $\cal{O}_{T}$ contains an immersed $(q,q,3)$ turnover, where $q \geq 4$.
	\item $T[2,2,4;n,3,r]$. $\cal{O}_{T}$ contains an immersed turnover of type $(2,4,r\geq5)$ (as well as the additional non-maximal turnovers listed in item (4)) if $n=2$, an immersed turnover of type $(4,4,r\geq3)$ if $n=3$, and immersed turnovers of types $(3,3,5)$, $(3,5,5)$ and $(5,5,5)$ if $n=2$ and $r=5$.
	\item $T[2,3,q;2,3,r]$. $\cal{O}_{T}$ contains an immersed $(q,r,r)$ turnover, where $q\geq3$ and $r=4$ or $r=5$.  
	\item $T[2,2,q;3,5,2]$. $\cal{O}_{T}$ contains an immersed $(q,q,5)$ turnover, where $q \geq 3$.
	\item $T[2,2,5;2,3,5]$. $\cal{O}_{T}$ contains an immersed $(3,5,5)$ turnover.
	\item $T[2,2,3;3,p,2]$. $\cal{O}_{T}$ contains immersed turnovers of type $(3,p,p)$ and $(p,p,p)$, where $p=5$ or $p=6$ (also, $(2,p,p)$ by item (5) and $(3,3,5)$, when $p=5$, by item (11)).
	\item $T[2,2,3;2,p,3]$. $\cal{O}_{T}$ contains immersed turnovers of type $(2,p,p)$, $(3,3,p)$ and $(p,p,p)$ if $p=5$, and an immersed turnover of type $(3,p,p)$ if $p=6$.
\end{enumerate}

\bibliographystyle{hyperamsplain}
\bibliography{refs}

\end{document}